\documentclass[11pt]{amsart}
\usepackage{graphicx,amssymb,amsfonts,amsmath,amsthm,newlfont}
\usepackage{epsfig,url}
\usepackage{color}
\usepackage{blkarray}
\usepackage{enumitem}
\usepackage{fancyhdr}

\usepackage[all,2cell]{xy} \UseAllTwocells \SilentMatrices

\newtheorem{thm}{Theorem}[section]
\newtheorem{cor}[thm]{Corollary}

\theoremstyle{definition}
\newtheorem{defn}[thm]{Definition}
\newtheorem{conj}{Conjecture}

\newtheorem{ex}[thm]{Examples}
\newtheorem{example}[thm]{Example}

\theoremstyle{remark}
\newtheorem{rem}[thm]{Remark}

\numberwithin{equation}{section}

\newcommand{\Z}{\mathbb Z}
\newcommand{\C}{\mathbb C}

\newcommand{\R}{\mathbb R}

\newcommand{\Pro}{\mathbb P}

\newcommand{\be}{\mathbf{e}}

\newcommand{\gr}{\mathrm{gr}}

\newcommand{\MT}{\mathcal{MT}}
\newcommand{\M}{\mathcal{M}}

\newcommand{\uu}{\mathfrak{u}}

\newcommand{\GG}{\mathcal{G} }
\newcommand{\UU}{\mathcal{U} }

\newcommand{\zetam}{\zeta^{ \mathfrak{m}}}

\newcommand{\Q}{\mathbb Q}

\newcommand{\F}{\mathbb F}

\newcommand{\To}{\longrightarrow}

\newcommand{\G}{\mathbb{G}}

\newcommand{\tone}{\overset{\rightarrow}{1}\!}

\newcommand{\dR}{\mathfrak{dr}}

\newcommand{\SL}{\mathrm{SL}}

\newcommand{\Or}{\mathcal{O}}

\newcommand{\mm}{\mathfrak{m} }

\newcommand{\id}{\mathrm{id} }

\newcommand{\per}{\mathrm{per}}

\newcommand{\dch}{\mathrm{dch}}

\begin{document}
\author{Francis Brown}

\begin{title}[From Deligne-Ihara  to multiple modular values.]{From the Deligne-Ihara conjecture to multiple modular values.}\end{title}

\maketitle

\section{Introduction}
I first met Professor Ihara at a conference in Kyoto on Galois-Teichm\"uller theory  in 2010. My talk was one of the very last of the conference
and had a particularly uninspiring title. I was  greatly honoured that Professor Ihara  patiently attended my talk, in which I  chaotically  sketched a proof of the conjecture attributed to him and Deligne.  It has been  a huge privilege, and a great pleasure,  to return to Kyoto in 2018 and participate in  the celebrations for his  80th birthday conference.

Whilst already mindful of the  profound impact of Ihara's ideas on my own work,  I was particularly struck at the conference by the scale of  his influence over  number theory as a whole, especially in Japan.  In this talk, I  shall only  focus on a tiny  fraction of Ihara's impressive legacy, namely his work on the projective line minus three points. It has two, related,  strands:
\begin{enumerate}
\vspace{0.1in}
\item (Genus $0$). I will report  on recent progress  related to the Deligne-Ihara conjecture and raise some new questions inspired by it      (\S\ref{sectP1},\ref{sectDIupdate}).
\vspace{0.1in}
\item  (Genus $1$). I will describe in \S\ref{sectgenus1} the main features of a nascent theory in genus one, and explain how it relates to  the Ihara-Takao relation, where modular forms make their first unexpected appearance. 
Paragraphs \S\S\ref{sectMMV}-\ref{sectStructuredouble} 
  are devoted to numerical examples  to illustrate the main features of this somewhat abstract theory and  make it more widely accessible. 
\end{enumerate}

In the sequel to this talk, which is logically independent from it, I explain how the periods associated to $(1)$ and $(2)$, namely multiple zeta values, and multiple modular values,  can be subsumed into a more general definition of multiple $L$-values.

\section{A brief history of the Deligne-Ihara conjecture}  \label{sectP1}

\subsection{The Deligne-Ihara conjecture (pro-$\ell$ version).} In 1984, Grothendieck \cite{Esquisse} proposed studying the profinite  completion of the fundamental group  
of the projective line minus three points $X= \Pro^1\backslash \{0,1,\infty\}$
with its outer action of $G_{\Q}= \mathrm{Gal}(\overline{\Q}/\Q)$. Shortly afterwards,  in a series of influential papers \cite{Ihara1, Ihara2, Ihara3, Ihara4}, Ihara initiated the study of the Galois action  on the pro-$\ell$ completion:
$$\rho_{\ell} :  G_{\Q} \To \mathrm{Out} \,  \pi_1^{(\ell)}(X,x)$$
where $\pi_1^{(\ell)}(X,x)$ is the inverse limit of all  finite quotients of $\pi_1(X(\C),x)$ whose  order  is a power of any prime $\ell$.   This action does not depend on the base point $x$. 
Ihara  used it to define, for every $\ell$, a   decreasing filtration on the absolute Galois group itself:
$$G_{\Q}=  I_{\ell}^0 G_{\Q}  \ \supseteq \  I_{\ell}^1 G_{\Q} \ \supseteq  \ I_{\ell}^2 G_{\Q} \ \supseteq \ \ldots\ . $$
Elements in the subgroup $I^k_{\ell}$ are  those which act as inner automorphisms on the $(k+1)^{\mathrm{th}}$-step in the lower central series filtration on $\pi_1^{(\ell)}(X,x)$. The first graded quotient in this filtration gives back  the $\ell$-adic cyclotomic character:
$$\gr^0_{I_{\ell}}  G_{\Q}  \cong  \Z_{\ell}^\times\ .$$
Ihara showed that $[I^p_{\ell}, I^q_{\ell}] \subset I^{p+q}_{\ell}$ and defined the associated graded  object
$$\mathfrak{g}^{\ell} = \bigoplus_{m \geq 1}  \gr^m_{I_{\ell}} G_{\Q}\ ,$$
which is a Lie algebra over $\Z_{\ell}$.  He proved that the degree $1$ part of this Lie algebra contains the images of Soul\'e characters 
$\sigma^{\ell}_{3}, \sigma^{\ell}_{5}, \ldots $
indexed by  every odd integer $\geq 3$, and made the following conjecture, jointly attributed to Deligne (\cite{Ihara3}, p. 600) see also  (\cite{Go2}  Conj. 2.1, \cite{Drinfeld}, p. 859):   
\begin{conj}  \label{conjDI} (Deligne-Ihara.) The Lie algebra $\mathfrak{g}^{\ell}\otimes_{\Z_{\ell}} \Q_{\ell}$ is freely generated  by the 
 elements $\sigma^{\ell}_{2n+1}$, for $n\geq 1$.
\end{conj} 

What makes this subject  so compelling is a  tension between the freeness in the previous conjecture, and the existence of infinitely many  `near-relations', the first of which occurs in weight $12$ and is  due to Ihara and Takao:
\begin{equation} \label{Ihara-Takao} 
\left[\sigma^{\ell}_3, \sigma^{\ell}_9\right]  -  3 \left[\sigma^{\ell}_5, \sigma^{\ell}_7\right]= \frac{691}{144} \, \delta^{\ell}_{12}\ .
\end{equation}
It turns out that  $\delta^{\ell}_{12}$, which can be represented as an element in a free  Lie algebra in two non-commuting elements,  is nearly zero: almost all  its  coefficients vanish (see \cite{YI}  and \cite{BrDepth}, \S8 for a precise statement).   This clearly makes life difficult if one wants to prove conjecture \ref{conjDI}, which implies that $\delta^{\ell}_{12}$ is non-zero.  Furthermore,  \eqref{Ihara-Takao} becomes an actual relation
in two different ways: firstly   modulo the prime 691, and  secondly  modulo  the depth filtration, since  $\delta^{\ell}_{12}$  vanishes, unexpectedly,  in  depths two (and  even three).  Therefore the modulo $\ell$ and depth-graded  analogues of conjecture \ref{conjDI} are  false: neither   $\mathfrak{g}^{\ell} \otimes_{\Z_{\ell}} \F_{\ell}$, nor $\gr_{\mathcal{D}} \left(\mathfrak{g}^{\ell}\otimes_{\Z_{\ell}} \Q_{\ell}\right)$  is freely generated  by the $\sigma_{2n+1}^{\ell}$. The element $\delta^{\ell}_{12}$ seems to be  related in a fundamental way to the existence of the first non-zero cusp form for $\SL_2(\Z)$ in weight 12.  Much of the work described in this talk  has been devoted to  gaining a better understanding of this  fascinating phenomenon.

\subsection{Motivic fundamental group}
A solution to the Deligne-Ihara conjecture came from a different part of mathematics: namely the theory of periods and multiple zeta values. The first step is to  formulate a motivic generalisation of  conjecture \ref{conjDI} and recast it in a different realisation. 

Based on Deligne's foundational work \cite{DePi1}  on realisations of the motivic fundamental group of $X=\Pro^1\backslash\{0,1,\infty\}$ and subsequent advances in the theory of motives, Deligne and Goncharov \cite{DeGo} defined a motivic fundamental groupoid:
$$\pi_1^{\mathfrak{m}} (X,\tone_0, - \tone_1)$$
as a pro-object in the abelian category of mixed Tate motives over $\Z$, denoted $\MT(\Z)$. 
The basepoints are tangent vectors of length $1$ (respectively $-1$) at the point $0$ (respectively $1$). One reason for these apparently obscure basepoints is that $X$ has no ordinary  points over $\Z$.  Another reason is that they formalise the logarithmic regularisation of divergent integrals. As a first approximation, the reader can pretend that they are the points $0$ and $1$ respectively, although these points lie in $\Pro^1$ and not in fact in $X$.

The motivic fundamental groupoid  admits the following classical realisations:
\begin{itemize}
\item (Betti). Its Betti realisation is the scheme over $\Q$ defined by the unipotent completion of the topological fundamental torsor of paths: 
$$\pi_1 (X,\tone_0, - \tone_1) \ . $$ It contains a distinguished path $\dch$ which is given by the straight line  from $0$ to $1$ along the real axis.   The $k$-points of the  Betti fundamental groupoid are given by the set of group-like formal power series 
$S \in k \langle\langle x_0,x_1\rangle \rangle $
in two non-commuting variables $x_0,x_1$, corresponding to loops around $0,1$. 
\vspace{0.05in}

\item (de Rham). The de Rham realisation is   dual to the $\Q$-vector space of iterated integrals on $\Pro^1\backslash\{0,1,\infty\}$. More precisely,
its affine ring  is isomorphic to the graded tensor coalgebra over $\Q$ on generators denoted by  
$$\omega_0 = \frac{dz}{z} \quad , \quad \omega_1 = \frac{dz}{1-z} \quad \in  \quad  \Gamma(\Pro^1, \Omega^1_{\Pro^1/\Q}(\log \{0,1,\infty\}))\ . $$
Its $k$-points are also isomorphic to the set of group-like formal power series $S \in k \langle\langle e_0,e_1\rangle \rangle $
in two non-commuting variables $e_0,e_1$ dual to  $\omega_0, \omega_1$.

\item  ($\ell$-adic).  The  set of $\Q_{\ell}$-points of its $\ell$-adic realisation  
$\pi_1^{\ell} (X,\tone_0, - \tone_1) (\Q_{\ell})$ is isomorphic to the group-like formal power series in two variables $\Q_{\ell} \langle\langle x_0, x_1 \rangle \rangle$ with a continuous action of $G_{\Q}$. There is a natural isomorphism  of $\Q_{\ell}$-schemes (\cite{HaMa}, \S A):
$$\pi_1^{\ell} (X,\tone_0, - \tone_1)  \cong\left( \pi^{(\ell)}_1(X,\tone_0, - \tone_1)\right)^{\mathrm{un}}$$
 where the right-hand side is the (continuous) unipotent completion of the pro-$\ell$ completion. The pro-$\ell$ completion of the topological fundamental groupoid  is  Zariski-dense in the $\Q_{\ell}$-points of the latter.  
\end{itemize} 

The Betti and de Rham realisations are related via the period isomorphism, and the theory of iterated integrals.  More precisely, the periods of the motivic fundamental groupoid include  the regularised iterated integrals  from $0$ to $1$:
$$\int_{\mathrm{dch}} \omega_{i_1} \ldots \omega_{i_n} \qquad \hbox{ for }\quad i_1,\ldots, i_n \in \{0,1\} \ . $$
By a remark due to Kontsevich, these can be expressed as $\Q$-linear combinations of multiple zeta values, which  are defined 
for $n_1,\ldots, n_{r-1}\geq 1$ and $n_r \geq 2$ by
\begin{equation} \label{MZVdef} \zeta(n_1,\ldots, n_r)= \sum_{1 \leq k_1 < \ldots <k_r}  \frac{1}{k_1^{n_1}\ldots k_r^{n_r}}
\end{equation}
and go back to Euler, at least for $r=2$. 
The $\ell$-adic theory, related to Betti via the Betti-$\ell$-adic comparison isomorphism,  
enables us to make the connection with conjecture \ref{conjDI} after taking $x$ to be one of the above tangential basepoints. There is also a crystalline theory which is related to the de Rham realisation via  $p$-adic periods.  The motivic structure of the fundamental groupoid is reflected in all of these realisations, but the theory of periods is perhaps the most accessible. For example, the Ihara-Takao relation is related to (but not in fact directly equivalent to)  the following elementary relation between double zeta values which first occurs in weight $12$ 
\begin{equation} 
28\, \zeta(3,9) + 150\, \zeta(5,7)+168\, \zeta(7,5) = \frac{5197}{691} \zeta(12) \ .\end{equation}
This relation is also directly related to the Ramanujan cusp form of weight $12$.
It could have been discovered by Euler, but  was found much more recently  \cite{GKZ}.

\subsection{Motivic reformulation  of  the Deligne-Ihara conjecture}
Since mixed Tate motives form a Tannakian category, the information contained in  the motivic fundamental groupoid is entirely encoded by any one of its realisations, together with the action of the corresponding Tannaka group. For example, 
$$\pi_1^{\bullet}( X, \tone_0, -\tone_1   ) \quad \hbox{ with its action of } \quad  
G^{\bullet}_{\MT(\Z)} = \mathrm{Aut}_{\MT(\Z)}^{\otimes}(\omega_{\bullet})  $$
where $\bullet \in \{B, dR, \ell\}$. It is known by the general theory, and most crucially Borel's theorems \cite{Borel1,Borel2} on the algebraic $K$-theory of $\Q$, that  $G^{\bullet}_{\MT(\Z)}$  is an extension
$$ 1 \To  U_{\MT(\Z)}^{\bullet} \To G^{\bullet}_{\MT(\Z)} \To \G_m \To 1$$
of the multiplicative group by a pro-unipotent group whose graded Lie algebra is the 
 free Lie algebra on non-canonical generators 
$\sigma^{\bullet}_3$,  $\sigma^{\bullet}_5, \ldots$
in every odd degree $-3,-5,\dots$. In the $\ell$-adic case,  there is a canonical continuous Zariski-dense homomorphism
$$G_{\Q} \To G^{\ell}_{\MT(\Z)}(\Q_{\ell})\ .$$
 It is easier to work in  the de Rham setting, where the above exact sequence   is canonically split.  
 The following theorem was conjectured by Goncharov and proved in \cite{MTZ}. 

\begin{thm} \label{thmDI} (Motivic version of the Deligne-Ihara conjecture).   The affine group scheme $G^{dR}_{\MT(\Z)}$ acts faithfully on $\pi_1^{dR}(X, \tone_0, -\tone_1)$.
\end{thm} 

Equivalently, the $\{\sigma^{dR}_{2n+1}\}$ act freely on   $\pi_1^{dR}( X, \tone_0, -\tone_1   )$. It follows that 
$$G^{\bullet}_{\MT(\Z)} \hbox{ acts faithfully  on the group scheme }   \pi_1^{\bullet}(X, \tone_0) $$
in any realisation $\bullet$. Taking $x=\tone_0$ and $\bullet = \ell$  proves conjecture \ref{conjDI}.

 \subsection{Classification of mixed Tate motives over $\Z$}
 The previous theorem  implies that every mixed Tate motive over $\Z$ is classified in the sense that it   appears inside the motivic fundamental groupoid of the projective line minus 3 points.
\begin{cor} The motivic fundamental groupoid $\pi_1^{\mathfrak{m}} (X,\tone_0, - \tone_1)$ generates the category $\MT(\Z)$. In other words, every mixed Tate motive over $\Z$ is isomorphic to a direct sum of  Tate twists of    subquotients of  its affine  ring.
\end{cor}

Such a classification theorem has applications. For instance:

\begin{cor} The periods of any mixed Tate motive over $\Z$ are $\Q[(2\pi i)^{-1}]$-linear combinations of multiple zeta values. 
\end{cor} 
There are many examples of integrals,  notably in high-energy physics, which are extremely hard to compute. When  the underlying geometry is mixed Tate over $\Z$,  one can deduce without computation that the integral is a multiple zeta value.

\subsection{Questions} 
In this talk, we   describe a precise relationship between cusp forms for $\SL_2(\Z)$ and the fundamental group of the projective line minus three points. 
This is directly  inspired by the   Ihara-Takao relation. 

In the sequel to this talk, we ask  whether 
 modular forms are intrinsically related to mixed Tate motives  at all, or whether this is 
 merely an artefact of the way  in which they are generated via $\Pro^1\backslash \{0,1,\infty\}$. 
More specifically, we consider whether  there exists a more  natural  way to construct the periods of mixed Tate motives without reference to $\Pro^1\backslash \{0,1,\infty\}$, which is seemingly  pulled  out of thin air.

\section{An update on the  Deligne-Ihara conjecture} \label{sectDIupdate}
Much of the combinatorial difficulty in understanding  the objects of the previous paragraph is that the  indexation is poorly adapted to the underlying motivic structure. I  explain how this can be rectified,  indicate some of the ingredients in the proof 
of  conjecture \ref{conjDI}, and  then turn to the depth-defect and Ihara-Takao relation. Some of the work described here  has occurred since  the paper \cite{MTZ} and  the surveys  \cite{DB, BrICM} appeared. 

\subsection{Decomposition of multiple zeta values}
One of the crucial  ingredients in the proof of  conjecture \ref{conjDI} is the notion of motivic period (or motivic multiple zeta value). We shall not spend much time motivating or defining this concept here, referring instead to \cite{BrICM}, \cite{BrNotes}. 
 In short, a motivic multiple zeta value is formally defined  for
 $ n_1,\ldots, n_{r-1} \geq 1$,  $n_r\geq 2$ to be an equivalence class
 $$\zetam(n_1,\ldots, n_r)  \ = \  [ \Or( \pi^{\mm}_1(\Pro^1\backslash \{0,1, \infty\},\tone_0, - \tone_1)), \omega_{\underline{n}}, \mathrm{dch}]^{\mm}$$
 which cuts out a `piece' of the affine ring of the fundamental groupoid. The form $\omega_{\underline{n}}$ is  the integrand of the expression for $\zeta(n_1,\ldots, n_r)$ as an integral along $\mathrm{dch}$.   These symbols satisfy the standard algebraic relations for multiple zeta values. They refine  an earlier notion due to Goncharov, for which the analogue of $\zeta(2)$ vanishes. 
 
 We only need to know that motivic multiple zeta values form a   $\Q$-algebra $\mathcal{H}$,  which is
 graded by the weight $n_1+\ldots+n_r\geq 0$, and   admits a period homomorphism:
\begin{eqnarray}  \mathrm{per}: \mathcal{H} &  \To&  \R \nonumber \\
\zetam(n_1,\ldots, n_r) & \mapsto& \zeta(n_1,\ldots, n_r) \ .\nonumber 
\end{eqnarray}  
It is conjectured to be injective. 
Motivic multiple zeta values come equipped with a certain `motivic' coaction, which is an  analogue of  Ihara's formula for the action of 
$\mathrm{Gal}(\overline{\Q}/\Q)$ on the pro-$\ell$ fundamental group. 
It takes the form 
\begin{eqnarray} 
\Delta: \mathcal{H}  & \To & \mathcal{H}   \otimes_{\Q}  \mathcal{A}\nonumber \\
\Delta ( \zetam(n_1,\ldots, n_r) ) & =&  \sum_I  c_I \,\zetam(a'_I) \otimes  \zeta^{\dR}(a_I)   \nonumber 
\end{eqnarray} 
where the sum is over some indexing set $I$, the $c_I \in \Q$, and $\zeta^{\dR}$ denote `de Rham' zeta values, which are defined in a similar manner to motivic multiple zeta values and can be thought of as `motivic multiple zeta's modulo $\zetam(2)$'. They generate an algebra $\mathcal{A}$   isomorphic to $\mathcal{H}/\zetam(2)\mathcal{H}$.   There exists a formula for this coaction, generalising an earlier one due to Goncharov, which \emph{a postiori} turns out to be exactly dual to a formula due to Ihara (see \cite{BrICM} for a quick derivation of the coaction formula from Ihara's).  One  can compute  the coaction (or parts of it) in practice \cite{MZVdecomp}.

As a  rough first approximation, one can think of a motivic multiple zeta value  as a \emph{matrix  of multiple zeta values} whose entries satisfy certain conditions, together  with a
distinguished entry.  The period homomorphism  picks out the distinguished entry; the coaction formula is  obtained by deleting rows and columns from the matrix.

 Using the coaction, together with some structure theorems on the category of mixed Tate motives,  we can decompose a motivic multiple zeta value into simpler pieces. 

\begin{thm}  Let $\Q\langle f_3,f_5,  \ldots \rangle$ denote the graded $\Q$-vector space generated by words in  non-commuting symbols $f_3,f_5, \ldots$ of weights $3,5,\ldots$  with the shuffle product.  Let $f_{2}$ denote a symbol of weight $2$ which commutes with the symbols $f_{2n+1}$ with  odd indices. 
There is a non-canonical injective homomorphism 
\begin{equation}\label{phimap}  \phi : \mathcal{H}  \To  \Q\langle f_3, f_5, \ldots \rangle \otimes_{\Q} \Q[f_2]
\end{equation}
which satisfies $ \phi(\zetam(2))  =  f_2$ and 
$$\phi( \zetam(2n+1) )   = f_{2n+1} \qquad \hbox{ for all } n  \geq 1\ . \nonumber $$
In particular, every motivic multiple zeta value admits, via $\phi$, a decomposition into the alphabet of letters $f_{2n+1}$ which respects all algebraic relations between the motivic multiple zeta values, and determines it uniquely. 
\end{thm}

As far as I know, an $f$-alphabet representation is presently the only clear way to navigate the complicated  relations satisfied by multiple zeta values.  This language will be used in an essential way when we discuss multiple modular values.

The  number of letters of odd weight   $f_{\mathrm{odd}}$ in a   word  defines an increasing filtration called the  \emph{coradical} or \emph{unipotency} filtration on $\mathcal{H}$ which we   denote by $C$.  It is well-defined (a conceptual and more general definition is given in \cite{BrNotes}, \S2.5).

\begin{example} The map $\phi$ depends on some choices. For some such choice,
\begin{eqnarray} \phi(\zetam(3,3)) &= & f_{3}f_{3}  - \frac{4}{35} f_2^3   \nonumber \\
\phi(\zetam(3,5))  &= &   - 5 f_3 f_5      \nonumber \\
\phi(\zetam(5,5)) & =& f_{5}f_{5} - \frac{16}{385} f_2^5  \nonumber \\
\phi(\zetam(7,3)) & =&  15 f_3f_7+6 f_5f_5+f_7f_3 - \frac{32}{385} f_2^5  \nonumber 
\end{eqnarray} 
For details on how to compute these formulas, see \cite{MZVdecomp}. The part of highest length in the symbols $f_{\mathrm{odd}}$ is canonically defined, i.e.,   the rational coefficient of the power of $f_2$ depends on the choice of  the homomorphism $\phi$, but not the rest. For example,
$$ \gr^C_2 (\phi) \, \left(  \zetam(7,3)\right) =   15 f_3f_7+6 f_5f_5+f_7f_3 $$
is independent of the choice of map $\phi$.   The interpretation  of this  formula is that $\zeta(7,3)$ is a period of a rank $5$ iterated extension of the following  Tate motives:
$$\Q \ , \  \Q(-3) \ , \  \Q(-5)\ , \  \Q(-7)\ , \ \Q(-10)\ , $$
i.e.,  $\zeta(7,3)$ is  an entry in a corresponding $5\times 5$ period matrix. The main  point of motivic periods is that they enable us to manipulate  such complicated  objects with great ease, avoiding  lengthy computations in homological algebra.   \end{example}

\subsection{The proof of the  Deligne-Ihara conjecture and recent developments}

\begin{thm}  \label{theoremphiisom} Any map  \eqref{phimap} is surjective. In other words, there is an isomorphism 
$$\phi: \mathcal{H}  \overset{\sim}{\To} \Q\langle f_3, f_5, \ldots \rangle \otimes_{\Q} \Q[f_2]\ .$$
It follows that every real, effective motivic period of a  mixed Tate motive over $\Z$ arises as a motivic multiple zeta value.  We deduce that every such motive occurs in the fundamental groupoid of $X$. Theorem \ref{thmDI} follows from Tannakian arguments. 
\end{thm}

Some comments:

\begin{itemize}
\item The proof of surjectivity  in \cite{MTZ}  goes by   showing that the Hoffman  elements
\begin{equation}\label{motivicHoffman} \zetam(n_1,\ldots, n_r) \qquad \hbox{ where } \quad n_i \in \{2,3\}
\end{equation}
have independent images under any $\phi$. They therefore form a basis for $\mathcal{H}$, which, by applying the period map $\mathrm{per}$, implies  a conjecture due to  Hoffman \cite{Hoffman}   stating that the $\zeta(n_1,\ldots, n_r)$ with $n_i \in \{2,3\}$  span the space of multiple zetas over $\Q$. The proof  uses in an essential way  a theorem   due to Zagier  \cite{Zagier23} who evaluated
$\zeta(2,\ldots, 2, 3, 2, \ldots,  2)$
as an explicit  linear form in $\zeta(2n+1)$ whose  coefficients are rational multiples of powers of $\pi$.  The coefficients in this linear form play an important role in the  proof. Zagier's original proof was a highly ingenious combination of combinatorial and analytic methods.
There have subsequently been several new proofs of this theorem (e.g., \cite{Zeta23v1},  \cite{Zeta23v2},  \cite{Zeta23v3}).  

In her thesis, Claire Glanois \cite{Glanois} constructed a new basis  for $\mathcal{H}$ using a clever variant of motivic Euler sums (iterated integrals on $\Pro^1\backslash \{0,-1,1,\infty\}$) which actually turn out to be motivic multiple zeta values.
\vspace{0.05 in}

\item The proof  is by induction on  a filtration by the `number of 3's' on the  elements \eqref{motivicHoffman}. This filtration, and the Hoffman basis itself, seems  mysterious at first sight, but is in fact naturally induced by a `block filtration'  on motivic multiple zeta values \cite{BlockFilt}  which is the filtration associated to the  block decomposition introduced by S. Charlton \cite{Charlton}. Adam Keilthy has  shown in his forthcoming Ph.\!\!\! D thesis that the block   filtration extends  to all motivic iterated integrals on $\Pro^1\backslash \{0,1,\infty\}$  and proved   many beautiful properties of it. 
\vspace{0.05in}

\item It is better to rephrase the previous theorem  after taking the associated  graded for  the coradical filtration. Thus theorem \ref{theoremphiisom} is equivalent to proving that  the  canonical  homomorphism $\gr^C \phi$ is  an  isomorphism: 
\begin{equation}\label{gradedphi} 
\gr^C (\phi)  : \gr^C \mathcal{H} \overset{\sim}{\To} \gr^C \Q\langle f_3, f_5, \ldots, \rangle \otimes_{\Q} \Q[f_2] \ .
\end{equation}
Stated in  this form, the proof can    be simplified, since several of the arguments actually hold for general motivic periods, as explained in \cite{BrNotes} . 
\vspace{0.05in}

\item The proof used  a mysterious $2$-adic argument  which can be  explained using  the $2$-adic Frobenius on the $\Z_2$-module generated by the motivic Hoffman basis. One  needs to apply the infinitesimal coaction to  elements
 \eqref{motivicHoffman} followed by   $\per_2 \otimes \id \pmod 2$, where $\per_2$ is the $2$-adic period homomorphism on de Rham zeta values. The key observation is that 
 $$v_2 \left(\zeta_2(3,2,\ldots, 2)\right) =1 $$
where $\zeta_2(3,2,\ldots, 2) = \per_2 \, \zeta^{\dR}(3,2,\ldots, 2) \in \Z_2$, and $v_2$ denotes the $2$-adic valuation.  This statement  follows from the calculation   $v_2 (\zeta_2(2k+1)) = v_2(\frac{4^k}{2k})$ which can be proved using a formula due to \"Unver \cite{Unver} \S5.11. This will be discussed elsewhere. 
\vspace{0.05in}
\item  Deligne has  analogous results for the projective line minus certain roots of unity \cite{DeRoots}.   See also Glanois \cite{Glanois} for further progress in this direction.
These are  situations in which there is no Ihara-Takao phenomenon. As a result, the motivic structure of the fundamental groupoid is much simpler than for the projective line minus three points. 
\end{itemize}

\subsection{Inverse problem} \label{sectInverseProb}
Given a mixed Tate motive over $\Z$ with prescribed weights, what are  its possible periods? 
To answer this question, one must compute the \emph{inverse} of  \eqref{gradedphi}. Equivalently, 
given a word $w$ in the $f_{2n+1}$, find a linear combination of motivic multiple zeta values whose image under   $\phi$  is $w$ to leading order.

\begin{defn} \label{defnzetaprescribed} Let $a_1,\ldots, a_r\geq 1$. Let  us denote by 
$$\zetam_{2a_1+1,\ldots, 2a_r+1}  \quad  \in  \quad C_r \mathcal{H}$$
 any motivic   multiple zeta value of weight   $2a_1+ \ldots +2a_r +r$
   whose leading term in any  $f$-alphabet decomposition \eqref{phimap}  is 
$ f_{2a_1+1}\ldots f_{2a_r+1}$. 
 In other words,
 $$\gr_r^C \phi ( \zetam_{2a_1+1,\ldots, 2a_r+1}) \  =  \  f_{2a_1+1}\ldots f_{2a_r+1} \ .$$
By theorem  \ref{theoremphiisom}  it exists.   It is only well-defined up to motivic multiple zeta values of coradical filtration $\leq r-1$ (in fact, one can do better: $\leq r-2$).   We shall write
$$ \zeta_{2a_1+1,\ldots, 2a_r+1} = \mathrm{per} \,  \left(\zetam_{2a_1+1,\ldots, 2a_r+1}\right)   $$
for some choice of element  $\zetam_{2a_1+1, \ldots, 2a_r+1}$. 
It is only well-defined modulo multiple zeta values  of   coradical filtration $\leq r-2$.  We have  $\zeta_{2n+1}= \zeta(2n+1)$. 
\end{defn}

\begin{ex} The usual indexation for multiple zeta values (and especially the Hoffmann multiple zeta values) are very ill-suited  to this problem, as a glance at the following examples shows. Here are some choices of representatives for $\zeta_{a,b}$:
 \begin{eqnarray} \label{zeta73etal}
 \zeta_{3,3} & = & \zeta(3,3)  \\
 \zeta_{5,3} & = & - \frac{1}{5} \zeta(3,5) \nonumber \\ 
  \zeta_{3,7} &=  & \zeta(3)\zeta(7)+  \frac{1}{14} \left(\zeta(3,7) + 3\, \zeta(5)^2\right) \nonumber 
 \end{eqnarray} 
 They are well-defined up to addition of a power of $\pi$. 
All hell breaks loose in weight twelve, since we are forced to include multiple zeta values of depth four: 
\begin{multline} \label{f93formula} \zeta_{3,9} \quad =  \quad  \frac{1}{19.691}  \Bigg( 2^4 3^2\,  \zeta(5, 3, 2, 2)  
- \frac{3^3  5.179   }{2.7}  \zeta(5, 7) -2.3^3 29\,  \zeta(7, 5)      \\
\qquad \qquad - \quad 3 .7^2 \zeta(3)\zeta(9) + 2^43\,  \zeta(3)^4 + 
 2^5 3^3 11 \,   \zeta(3, 7) \zeta(2)  +   2^5 3^2 31 \, \zeta(7, 3)  \zeta(2)   \\
\quad  \qquad  -  \quad  2^4 3^4 \, \zeta(3, 5) \zeta(4)  -   2^5 3^2\,  \zeta(5, 3)  \zeta(4) - 2^3 3. 5^2 \zeta(3)^2 \zeta(6)  + \frac{3.128583229}{2^4 7. 691} \, \zeta(12) 
\Bigg)
\end{multline}
This number  is  only well-defined  modulo $ \pi^{12} \Q$ (but this particular representative will appear as a multiple modular value later).  The  problem `write down the periods of all mixed Tate motives over $\Z$ with three weight-graded pieces $\Q, \Q(-3), \Q(-12)$'  is  therefore surprisingly complicated. A period matrix of such an object  can be written in the form
$$
\begin{pmatrix}
1      & \zeta_3 & \zeta_{3,9} \\
 0 & (2\pi i)^3 & (2\pi i )^3 \zeta_9 \\
 0 &  0   & (2\pi i)^{12} \\
\end{pmatrix}
$$
with respect to suitable bases. This  illustrates the combinatorial complexity which implicity hides behind  theorem \ref{thmDI}, and is due to the Ihara-Takao phenomenon.
\end{ex}

\subsection{Ihara-Takao and the depth-defect}
The depth filtration on motivic multiple zeta values is compatible with the coradical filtration:
$$\zetam(n_1,\ldots, n_r) \quad  \in  \quad  C_r \mathcal{H}\ .$$
Thus, an $f$-alphabet decomposition \eqref{phimap} of a motivic multiple zeta value of depth $r$ involves words in the letters $f_{2n+1}$
of length at most $r$. Alas, the converse is false: the  depth filtration is strictly smaller than the coradical filtration. 
This phenomenon first becomes visible in weight $12$, where it is related to the fact that $\delta^{\ell}_{12}$ in \eqref{Ihara-Takao}
  vanishes in depth 2. 
For example, one can show that $\zeta_{3,9}$ and $\zeta_{5,7}$ cannot individually be expressed as multiple zeta values of depth $2$. Only the following linear combination lies in the subspace generated by  double zeta values:
\begin{equation}  \label{f75andf93} 
 \zeta_{5,7}+ 3\, \zeta_{3,9} \quad =  \quad \frac{1}{9} \zeta(3, 9)+3 \,  \zeta(3) \zeta(9)+\frac{5}{3}\zeta(5)\zeta(7)-\frac{31. 139}{2.691}\zeta(12)  \ .
 \end{equation} 
   The expression on the left-hand side is  dual to the left-hand side of \eqref{Ihara-Takao}.

\section{A sketch of the theory in genus one} \label{sectgenus1} 
The element  $\delta^{\ell}_{12}$ in the Ihara-Takao relation \eqref{Ihara-Takao}
occurs in the same weight as  the first cusp form for $\SL_2(\Z)$. This suggests  studying the fundamental group of the moduli stack $\mathcal{M}_{1,1}$ of elliptic curves, which, together with   $\mathcal{M}_{0,4} \cong \Pro^1\backslash \{0,1,\infty\}$,  comprise the two  `generators'  in Grothendieck's  tower  of moduli spaces  $\mathcal{M}_{g,n}$ of curves of genus $g$ with $n$ marked points \cite{Esquisse}.

The analytic space associated to  $\mathcal{M}_{1,1}$ is the orbifold quotient
$$\mathcal{M}_{1,1}(\C) = \Gamma \backslash\! \backslash \mathfrak{H}$$
where $\mathfrak{H}$ is the upper half plane $\mathfrak{H}=\{\tau \in \C: \mathrm{Im}\, \tau>0\}$  and $\Gamma= \SL_2(\Z)$.  A  natural choice for a basepoint is the unit tangent vector at the cusp, denoted by 
\begin{equation} \label{TBpoints} \partial/\partial q   \quad \hbox{ or } \quad \tone_{\infty}
\end{equation}
where $q=\exp( 2\pi i \tau)$. The former denotes the unit tangent vector on the punctured $q$-disc and induces a tangential basepoint on $\mathcal{M}_{1,1}$.  We tend to use  the latter notation to denote a tangential basepoint on the universal covering, i.e.,  a tangent vector of length $1$ on the tangent space at $\tau = i \infty$ on the compactified  upper half plane. They coincide on  $\mathcal{M}_{1,1}(\C)$.    This basepoint has an infinite $\Gamma$-stabiliser. 
Other natural choices of basepoint include the images of the points $\tau =i$ or $\tau = e^{2\pi i/3}$  whose $\Gamma$-stabilisers have order $4$ and $6$ respectively.   Better still, one can consider the groupoid of all three basepoints simultaneously and the paths between them. 

 There is a canonical isomorphism
$$\pi_1(\mathcal{M}_{1,1}(\C); \tone_{\infty}) = \Gamma \ .$$
The first problem one encounters is that, although the profinite completion of $\SL_2(\Z)$ is huge (e.g., by Belyi's theorem \cite{Belyi}),  its   unipotent completion  is trivial:
$$\pi_1^{\mathrm{un}}  (\mathcal{M}_{1,1}(\C); \tone_{\infty}) =1$$
There are many ways to see this: perhaps the simplest is that $H^1(\mathcal{M}_{1,1};\Q)$ vanishes as there are no modular forms of full level in weight $2$.
The remedy is to consider not the unipotent completion but a generalisation called the  \emph{relative completion} \cite{HaMHS} of $\Gamma$ with respect to the inclusion $\Gamma \subset \SL_2(\Q)$. This is a pro-algebraic group  over $\Q$
$$\GG^B_{1,1} = \pi_1^{\mathrm{rel}}(\mathcal{M}_{1,1}; \tone_{\infty}) $$ which is an extension  of the algebraic group $\SL_2$
$$1 \To \UU^B_{1,1} \To \GG^B_{1,1} \To \SL_2\To 1$$
 by a  pro-unipotent algebraic group $\UU^B_{1,1}$. It admits a Zariski-dense homomorphism $\Gamma \rightarrow \GG^B_{1,1}(\Q)$.  It can be thought of as an algebraic envelope of $\Gamma$ that is neither too small (as in the case of the unipotent completion) nor too large (as in the case of the pro-algebraic fundamental group), but lies in a Goldilocks zone  in-between.

\begin{itemize}
\item  The relative completion has various realisations $\GG_{1,1}^{\bullet}$, $\UU_{1,1}^{\bullet}$ which have a natural Tannakian description. In each case, one considers a certain Tannakian category of sheaves on $\mathcal{M}_{1,1}$ (e.g., local systems, or algebraic vector  bundles with integrable connection) which have an increasing filtration with  graded pieces  of a  specified type. One  defines  $\GG^{\bullet}_{1,1}$ to be the group of tensor automorphisms of the  functor `fiber at \eqref{TBpoints}'.  Most of the properties of $\GG_{1,1}^{\bullet}$ follow automatically. 
Since unipotent completion is a special case of relative completion, this construction also  applies  in the case of the projective line minus three points.
\vspace{0.05in}

\item  The de Rham description is particularly accessible. One shows \cite{HaGPS} that $\uu_{1,1}^{dR}  = \mathrm{Lie}\, \UU_{1,1}^{dR}$ is isomorphic to the completed free Lie algebra generated by
$$\prod_{n\geq 0} H^1_{dR}(\mathcal{M}_{1,1};  \mathbb{V}^{dR}_n)^{\vee} \otimes   V_n^{dR} $$
where $\mathbb{V}^{dR}_n= \mathrm{Sym}^n \mathcal{H}$, and
$\mathcal{H}=R^1 \pi_* \Omega^{\bullet}_{\mathcal{E}/\mathcal{M}_{1,1}}$  where $\pi:\mathcal{E}\rightarrow \mathcal{M}_{1,1}$  is the 
 universal  elliptic curve.
 It is an algebraic vector bundle  with the Gauss-Manin connection.  Here $V_n^{dR}$ denotes the fiber of $\mathbb{V}_n^{dR}$ at the basepoint \eqref{TBpoints}.
  In the theory of modular forms, one often identifies 
 $$V_{n}^{dR}= \bigoplus_{i+j=n}  X^i Y^j \Q $$
 with the space of homogeneous polynomials of degree $n$ in two variables $X$ and $Y$, with its natural right $\SL_2$-action.   By the Eichler-Shimura theorem, 
 $$\dim_{\Q} \, H^1_{dR}(\mathcal{M}_{1,1};  \mathbb{V}^{dR}_n) = 1 + 2 \dim_{\C} S_{n+2}(\Gamma)$$
 where $S_{n+2}(\Gamma)$ denotes the space of cusp forms for $\Gamma$.   In fact, there is a canonical  injection 
 $M_{n+2}(\Gamma;\Q) \rightarrow  H^1_{dR}(\mathcal{M}_{1,1};  \mathbb{V}^{dR}_n)$
 from the space of modular forms of level $1$ and weight $n+2$ with rational Fourier coefficients into this space. The remaining classes can be generated, for example, using weakly holomorphic modular forms (i.e., with a pole at the cusp).  It follows that  the Lie algebra $\uu_{1,1}^{dR}$ is isomorphic to the completed  free Lie algebra on symbols
 $$\be_{n+2} \otimes V^{dR}_{n}  \quad  , \quad  \be'_f \otimes V^{dR}_{n} \quad , \quad  \be''_f \otimes V^{dR}_{n}$$
  where $\be_{n+2}$ corresponds to an  Eisenstein series, and $\be'_f, \be''_f$ correspond to a cusp form (resp. weakly holomorphic cusp form)  of weight $n+2$. 
   The reader is warned  that unlike the case of the projective line minus three points, the de Rham relative completion is not canonically graded and the choices of generators are non-canonical.  The periods of relative completion (called `multiple modular values' in \cite{MMV}) include the regularised iterated integrals of   vector valued modular forms along the imaginary axis.  In the simplest case, these reduce to classical  Eichler integrals of cusp forms (see below). 
   
\vspace{0.05in}
\item  A key ingredient in the study of $\Pro^1\backslash \{0,1,\infty\}$ is a bound in the extension groups in the category of mixed Tate motives over $\Z$, which follows from Borel's deep results in algebraic $K$ theory. Such results are lacking in the modular situation, but we can work instead in a category of realisations (as in \cite{DePi1})  over $\Q$ with coefficients in $\overline{\Q}$. Hereafter we use the word `motive' loosely to mean an object in the subcategory  generated by  $\Or(\GG_{1,1})$. 
This formalism is sufficient to exhibit elements of a `motivic' Lie algebra generalising the algebra generated by the  $\sigma_{2n+1}$ which play a role in the Deligne-Ihara conjecture. Such elements (which are only defined up to commutators) correspond to non-trivial simple extensions of $\Q$ inside $\Or(\GG_{1,1})$. 
In contrast to the case of genus zero, we now expect infinitely many different types of such extensions.
\end{itemize} 

We briefly describe some of the `motivic' elements which act non-trivially on the relative completion $\GG_{1,1}$, and in particular the Lie algebra $\uu_{1,1}$.
\subsubsection{Zeta elements} 
 Analogues of the familiar  zeta elements  $\sigma_{2n+1}$ indeed appear, and correspond to extensions of Tate objects occuring inside $\Or(\mathcal{G}_{1,1})$: 
$$0 \To \Q(2n+1) \To \mathcal{E} \To \Q \To 0\ .$$
One can prove (\cite{MMV}, \S19) that  their de Rham versions act on $\uu^{dR}_{1,1}$ via
\begin{equation} \label{zetaelements} \sigma_{2n+1}^{dR} = \frac{(2n)!}{2} \,   \mathrm{ad} \left(  \be_{2n+2} Y^{2n}\right)  + \ldots 
\end{equation} 
In fact, some of the higher order terms are also known. This already has applications to the theory of the projective line minus $3$ points, as it enables us to construct zeta elements $\sigma_{2n+1}$ which are canonical up to depth four \cite{BrSigma}. 

\subsubsection{Modular elements}
These have no analogue in the case of the projective line minus three points. 
For any eigencusp form $f$ of weight $w$ and any integer $d\geq w$, there exist \emph{modular elements}  $\sigma_{f}(d)$ which correspond to a non-trivial  extension 
$$0 \To M_f(d) \To \mathcal{E} \To \Q \To 0  \ , $$
where $M_f$ is the simple object of rank $2$ (`motive' in the category of realisations) associated to the cusp form $f$. One proves \cite{MMV} that they act on $\uu^{dR}_{1,1}$    via
$$\sigma^{dR}_{f(d)} = \mathrm{ad} ( \mathbf{b}_{f(d)} ) + \delta_{f(d)}$$
where the `geometric part' $\mathbf{b}_{f(d)} \in [\uu^{dR}_{1,1},\uu^{dR}_{1,1}]$ 
and  the `arithmetic part'   $\delta_{f(d)} $ is an outer derivation of $ \uu^{dR}_{1,1}$.  Their leading terms  are known.
In the case $d=w$, we have 
\begin{equation} \label{modulardelta} 
\delta_{f(w)} ( \be'_f  Y^{w-2} )   =  \sum_{2a+2b= w-2} \,    c_f^{a,b} [ \be_{2a+2} Y^{2a}, \be_{2b+2} Y^{2b} ]  +\ldots 
\end{equation} 
where $c_f^{a,b}$ are algebraic numbers in the field $K_f$ generated by the Fourier coefficients of $f$, which are   proportional to the  coefficients in the even period polynomial of $f$. More precisely,    one has  the equality of  points  in projective space $\Pro^{n-1}(K_f)$
$$( c_f^{1,n-1} \ : \     \ldots \  : \  c_f^{n-1,1}   )  \ =  \ ( i^3 \Lambda(f,3) \  : \   i^5 \Lambda(f,5) \  : \  \ldots \  : \   i^{2n-1} \Lambda(f,2n-1) ) $$
where $w=2n+2$ and $\Lambda$ is the completed $L$-function of $f$. 

\subsubsection{Higher elements}
We expect $\sigma^{dR}_{2n+1}, \sigma^{dR}_{f(d)}$, \ldots  to be the beginning of an infinite sequence of  families of elements.  The next  family is $\sigma^{dR}_{f\otimes g(d)}$ corresponding to an extension of $\Q$ by a Rankin-Selberg motive $M_f \otimes M_g(d)$, for $d$ first in the `semi-critical' range, and later for all values of $d$. After this come extensions of $\Q$ by  $M_{f}\otimes M_g \otimes M_h(d)$,  and so on.

\subsection{A genus one Deligne-Ihara conjecture?} Beilinson's conjecture predicts exactly how many motivic extensions of $\Q$ by 
\begin{equation} \label{SymMmotive} \left(\mathrm{Sym}^{m_1} M_{f_1} \otimes \ldots \otimes \mathrm{Sym}^{m_r} M_{f_r}\right)(d) 
\end{equation} 
one should  expect to find in nature. Such extensions correspond to derivations
$$\sigma_{2n+1}^{dR} \  , \  \sigma_{f(d)}^{dR} \ , \ \sigma^{dR}_{f\otimes g(d)} \ , \ \sigma^{dR}_{(\mathrm{Sym}^2 f)(d)} \ , \   \ldots $$ 
which will be non-trivial if the extension actually occurs in $\Or(\GG_{1,1})$. One hopes  this is true, and furthermore, that  these derivations generate a free Lie algebra, by analogy with conjecture \ref{conjDI}.
 This  would    imply a classification theorem for mixed modular motives, assuming an analogue of Borel's theorem which would bound the extension groups of mixed modular motives.   As a first approximation, we can prove \cite{MMV}:

\begin{thm} Any representative for the  zeta and modular elements $\sigma^{dR}_{2n+1}, \sigma^{dR}_{f(d)}$ act freely upon $\uu^{dR}_{1,1}$, i.e., generate a free Lie algebra.  
\end{thm} 
In fact, one can prove much more, namely that any non-zero `motivic' derivations satisfying quite a weak condition necessarily generate a free Lie algebra \cite{MMV}, \S21.

Here, however, comes a surprise.   One can show that not every predicted extension of $\Q$ by  \eqref{SymMmotive}  can  actually arise in $\Or(\GG_{1,1})$: there are certain  exceptional cases for large $m_1+\ldots+m_r$ and small $d$  which \emph{cannot} occur due to a Hodge-theoretic obstruction.  This is very mysterious.

\begin{table}[h!]
 \renewcommand*{\arraystretch}{1.6}
\begin{center}
\begin{tabular}{|c|c|c|} \hline
Object  &  $\mathcal{M}_{0,4}= \Pro^1\backslash \{0,1,\infty\} $& $\mathcal{M}_{1,1}$\\ \hline
$\pi^{?}_1$  &  Unipotent completion  &  Relative completion  \\
  & $\pi_1^{\mathrm{un}}(\Pro^1\backslash \{0,1,\infty\})$ & $\GG_{1,1}=\pi_1^{\mathrm{rel}}(\mathcal{M}_{1,1})$ \\ \hline
 basepoints &  $\tone_0, -\tone_1$  &   $ \frac{\partial}{\partial q} $  \quad ($\tone_{\infty}$) \\  \hline
 Lie algebra of &  Free Lie algebra generated by    &    Free Lie algebra generated by \\ 
 $\pi^{dR}_1$  & $ e_0, e_1  $ dual to &  $\{\be_{w},  \be'_f ,  \be''_f \}   \otimes V^{dR}_{w-2}$   \\
  &   $\frac{dz}{z}, \frac{dz}{1-z}$    & for $f \in S_w(\Q)$ cusp form  \\
  \hline
  Paths &  The straight line path  & The images of the paths \\
  &    $\mathrm{dch}$  &   $S=\begin{pmatrix} 0 & -1 \\ 1 & 0 \end{pmatrix}$ and $T =\begin{pmatrix} 1 & 1 \\ 0 & 1 \end{pmatrix} $ \\
  &  from $\tone_0$ to $-\tone_1$ & i.e., $i \R_{>0}$ and a loop in $q$-disk \\
   \hline
  Periods &  Regularised iterated integrals of  & Regularised  iterated integrals of \\
  & one forms & vector-valued modular forms \\
  &  $\frac{dz}{z}$ and $\frac{dz}{1-z}$  &  of first and second kinds  \\
  & a.k.a. multiple zeta values &  a.k.a.  multiple modular values  \\ \hline
  Pure periods  & Polynomials in $2\pi i$ & Polynomials in $\omega_{f}^{+}, i \omega_f^{-}, \eta_f^{+}, i \eta_f^-$ \\  \hline
   `Motives' & Mixed Tate Motives over $\Z$ & Mixed Modular `Motives' \\ 
  &  $\MT(\Z)$ &     $\mathcal{MMM}_{\SL_2(\Z)}$ \\
 Simple objects  & $\Q(n)$  &    $\left(\mathrm{Sym}^{m_1} M_{f_1} \otimes \ldots \otimes \mathrm{Sym}^{m_r} M_{f_r}\right)(d) $ \\ \hline
 Motivic Lie  & \underline{Generators}     & \underline{Generators} \\
 algebra & $\zeta$-elements:  $\sigma_3,\sigma_5,\ldots $  &  $\zeta$-elements: $ \sigma_{3}, \sigma_5, \ldots $ \\
 & & Modular elements: $\sigma_{f}(d)$ \\
 &  & Rankin-Selberg elements, $ \cdots$ \\ \hline
  \end{tabular} 
\end{center}
\label{default}
\end{table}%

\subsection{The Ihara-Takao relation revisited}  Let $\Delta$ denote the Ramanujan cusp form of weight $12$.    The  critical values of its completed $L$-function satisfy 
$$( i^3 \Lambda(\Delta,3) : i^5 \Lambda(\Delta,5):  i^7 \Lambda(\Delta,7):  i^9 \Lambda(\Delta, 9)) =(14: -9: - 9:14) \ . $$
It follows that  some normalisation of the associated   modular  element  (which is only defined up to scalar multiple) $\sigma^{dR}_{\Delta(12)}$ acts via the derivation
$$  \be'_{\Delta} Y^{10} \quad \mapsto \quad   14\,  [ \be_{4} Y^2, \be_{10} Y^8] -  9\,  [ \be_6 Y^4, \be_{8} Y^6 ]   \  + \  \left(\hbox{higher order terms}\right)  $$
Therefore if we rescale this particular choice of generator and define 
$$
d^{dR}_{12} = 1440\,  \sigma^{dR}_{\Delta(12)} \left( \mathrm{ad} \, \be'_{\Delta} Y^{10} \right)$$
and substitute in  equation \eqref{zetaelements}, we obtain  the congruence
\begin{equation} \label{ModularIT} 
  d^{dR}_{12} \equiv  \left[ \sigma^{dR}_3,\sigma^{dR}_9 \right] - 3 \left[\sigma^{dR}_5, \sigma^{dR}_7\right]
  \end{equation} 
modulo higher order Lie brackets.  Compare with  \eqref{Ihara-Takao}. Note, however, that the previous equation takes place in a different Lie algebra, namely $\mathrm{Der} \,  \uu_{1,1}^{dR}$. To make the connection with the projective line minus three points, one can use a `monodromy' homomorphism from $\mathrm{Der}\, \uu_{1,1}^{dR}$  to the derivations on the de Rham $\pi_1$  of the  punctured infinitesimal Tate curve, which sends $d^{dR}_{12}$ to zero. Hain has shown how to  relate the latter to the projective line minus three points. We  can deduce, for example,  that 
$$ 0 \equiv [\sigma^{dR}_3,\sigma^{dR}_9] - 3 [\sigma^{dR}_5, \sigma^{dR}_7]$$
in the associated depth-graded Lie algebra \cite{MMV}.

A version of the above should hold $\ell$-adically as well. One expects  a `Galois' Lie algebra constructed out of a filtration on $G_{\Q}$ induced  by its action on the  lower central series of the $\ell$-adic realisation  $\GG^{\ell}_{1,1}$. One expects that $\ell$-adic versions of the modular elements $\sigma_{f(d)}^{\ell}$ should be non-zero. Repeating the above argument in the $\ell$-adic setting should explain the origin of $\delta^{\ell}_{12}$ in the Ihara-Takao relation $(\ref{Ihara-Takao})$. The fact that it vanishes modulo $691$ should be related to the congruence  $\Delta \equiv \G_{12} \pmod{691}$ and its $\ell$-adic interpretation via the action of $G_{\Q}$ on $(M_{\Delta})_{\ell}$.   This  seems to be related to the fascinating talk of Sharifi in this  volume.

\section{Multiple modular values of level $1$: an overview} \label{sectMMV}
The rest of this talk concerns multiple modular values, which are the periods of the relative completion of the fundamental group of $\mathcal{M}_{1,1}$ and the genus one analogues of multiple zeta values.  The totally holomorphic multiple modular values are very concrete and provide  an elementary perspective on the underlying motivic structure. 
After giving a brief overview of general principles, the remaining sections  consist of  examples which  may  help to build intuition  for  the  theory.

Let $f_1,\ldots, f_r$ be modular forms for the full modular group  $\SL_2(\Z)$ which have  rational Fourier coefficients.  Let $w_1,\ldots, w_r$ denote their weights. 

\begin{defn} A  totally holomorphic multiple modular value  of length $r$   is  $(2\pi i)^{w_1+\ldots +w_r-r}$ times a regularised iterated  integral   (\cite{MMV}, \S5)
\begin{equation} \label{Lambdaf1fr} \Lambda(f_1, \ldots,f_r ; n_1 ,\ldots, n_r) =     \int_0^{\tone_{\infty}}  
f_1(t) t^{n_1-1}dt  \ldots  f_r(t) t^{n_r-1} dt 
\end{equation}
where each $n_i$, for $1\leq i \leq r$, satisfies $0 < n_i < w_i$.  
\end{defn}

The case where all $f_i$ are cuspidal was considered by Manin, and the integral is  an ordinary integral from $0$ to $\infty$. In general, the tangential basepoint at the upper range of integration implies a certain regularisation procedure  described in the sequel to this talk. The numbers \eqref{Lambdaf1fr} are  periods of the relative completion of the fundamental group of $\mathcal{M}_{1,1}$, but not the only periods of the latter (see \S \ref{sectQP}).
\subsection{Properties} 

It follows from  the general properties of iterated integrals  that the  numbers \eqref{Lambdaf1fr} satisfy a reflection symmetry  and shuffle product relations. In particular,  they generate an algebra. In fact, much more is known:

\begin{itemize}
 \item The  values $\Lambda(f_1,\ldots, f_r;n_1,\ldots, n_r)$ for fixed $f_i$ and varying $0 < n_i < w_i$ satisfy non-abelian cocycle relations  (\cite{Ma1}, \cite{MMV} \S5).
 \vspace{0.05in} 
   \item Surprisingly, there is a mechanism of transference of periods  between integrals of \emph{different} sets of modular forms. It is a higher analogue of the Peterson inner product.  It generates relations between multiple modular values $$\Lambda(f_1,\ldots, f_i;   - ) \Lambda(f_{i+1}, \ldots, f_n; -) \qquad \hbox{ for } \quad  1\leq i < n$$
where $f_1,\ldots, f_n$ are any fixed set of  modular forms.
\vspace{0.05in}

 \item Certain linear combinations of iterated Eisenstein integrals are multiple zeta values. This is a consequence of the fact that the relative completion of the fundamental group of $\mathcal{M}_{1,1}$  has a monodromy representation on the pro-unipotent fundamental group of the punctured Tate curve, whose periods are multiple zeta values. Enriquez has written down a combinatorial relation between the Drinfeld associator and the periods of the latter \cite{Enriquez}. This is closely related to the  elliptic MZV's  developed by him and Matthes \cite{Enriquez2, Matthes}.
 \vspace{0.05in}
   \end{itemize} 
 This list is not  exhaustive:
Beilinson's conjectures on extensions of motives should imply  relations  which presently have no elementary or direct proof.

\begin{rem} The transference principle can be viewed as  a  manifestation of some kind of   convolution on modular forms.  For example, it implies a relation between 
multiple modular values of the form $\Lambda(f_1,f_2)\Lambda(f_3) $ and those of the form 
$\Lambda(f_1)\Lambda(f_2,f_3)$.
It follows that the multiple modular values $\Lambda(f_1,f_2;n_1,n_2)$ for  any two modular forms $f_1,f_2$ pick up information from \emph{a priori}   unrelated  modular forms $f_3$. This phenomenon is related in simple cases  to the Rankin-Selberg method, but is more general. 
\end{rem} 

\subsection{Two relations to multiple zeta values} We
describe two different connections with  multiple zeta values. One is conjectural, the other is a theorem. 

Recall that the Eisenstein series of weight $2k\geq 4$ is the modular form of level one and  weight $2k$ with Fourier expansion normalised as follows:
$$\G_{2k} = -\frac{b_{2k}}{4k} + \sum_{n\geq 1}  \sigma_{2k-1}(n) q^n\ ,$$
where $\sigma$ denotes the divisor function. Recall definition  \ref{defnzetaprescribed}.

 \begin{conj}  \label{conjEisasMZV}  Let $a_1,\ldots, a_r\geq 1$ and $a= a_1+\ldots +a_r$. Then 
$$(2\pi)^{2a+r} \Lambda(\G_{2a_1+2}, \ldots, \G_{2a_r+2}; 1,\ldots, 1)   \equiv   (-1)^a \, \frac{(2a_1\!+\!1)! \ldots (2a_r\!+\!1)!}{2^r} \,  \zeta_{2a_1+1,\ldots, 2a_r+1} \\
  $$
  modulo multiple modular values  of unipotency (coradical) filtration $\leq r-1$.  \end{conj}

We expect the following   stronger  variant  to be true:  the difference between the left and right-hand sides of the equation  are   linear combinations of multiple modular values of the form
  $\Lambda(\G_{2a'_1+2},\ldots, \G_{2a'_k+2}; n_1,\ldots, n_k)$ where $k\leq r$. In all the examples we can actually take $k=r$ and  $a'_1=a_1,\ldots, a'_k=a_k$. In other words,  we find that a representative for the multiple zeta value $\zeta_{2a_1+1,\ldots, 2a_r+1}$ can be expressed as a $\Q$-linear combination of  the multiple modular values  $\Lambda(\G_{2a_1+2}, \ldots, \G_{2a_r+2};n_1,\ldots, n_r)$.
     
      The motivation for this conjecture is theorem 22.2 in \cite{MMV}, where a similar statement is proved on the level of mixed Hodge structures. Note that iterated integrals of the Eisenstein series $\G_{2k}$   are \emph{not} periods of mixed Tate motives, and are not expected to be multiple zeta values in general (see the  examples below).
           So if conjecture \ref{conjEisasMZV} is true, then we have a solution to the inverse problem \S \ref{sectInverseProb}. The price to pay is the  possible introduction of non multiple-zeta value periods.

Now compare this with the following theorem, which  is a  corollary of a much more general and precise  theorem about motivic periods which will be part of A. Saad's forthcoming  doctoral thesis at the University of Oxford. 
\begin{thm}\label{thmSaad} (Saad) Every multiple zeta value of depth $r$ and weight $2a+r$ can be expressed  as a rational linear combination of iterated integrals of Eisenstein series 
$$(2\pi )^{2a+r} \Lambda(\G_{2a_1+2}, \ldots, \G_{2a_r+2}; n_1,\ldots, n_r) \quad \hbox{ where } \quad  1\leq n_i \leq 2a_i+1 \ .$$
\end{thm}

The relationship between this theorem and conjecture \ref{conjEisasMZV} is subtle: already in length $r=2$, conjecture \ref{conjEisasMZV}
relates  $\Lambda(\G_{4}, \G_{10};1,1)$ to the  number $\zeta_{3,9}$, which cannot be expressed as multiple zeta values of depth $\leq 2$ (we expect).   Indeed, it is not captured by the previous theorem for $r=2$ except in the linear  combination  \eqref{f75andf93}, but will reappear in depth $r=4$.

\section{Examples of multiple modular values in length one.} \label{sectMMVex1}
Multiple modular values in length one are generated by:   powers of $2\pi i$,  odd zeta values, periods of cusp  forms (all of which are totally holomorphic) and quasi-periods of cusp forms (which are not totally holomorphic). 

For $f $ a modular form of weight $2k$, the totally holomorphic values \eqref{Lambdaf1fr} are classical.  As the notation suggests, they are critical values of the completed $L$-function
 $$\Lambda(f; s) = (2\pi)^{-s} \Gamma(s) L(f,  s)  $$
where $L(f, s) = \sum_{n\geq 1} a_n(f) n^{-s}$, and $a_n(f)$ are the Fourier coefficients of $f$. The functional equation $\Lambda(f ; s) = (-1)^k \Lambda(f; 2k-s)$  holds.

\subsection{Eisenstein series.} 
Since the $L$-function of $\G_{2k}$ factorises as a product
$$L(\G_{2k}, s) = \zeta(s) \zeta(s-2k+1)$$
 of two Riemann zeta functions,  we deduce that  the  even values satisfy   $$\Lambda(\G_{2k};2i) =  \frac{(-1)^i}{2} \frac{b_{2i}}{2i} \frac{b_{2k-2i}}{2k-2i} \ .  $$ 
The values $\Lambda(\G_{2k}; n)$ vanish for odd values of $3 \leq n \leq 2k-3$  except for
 $$ (2 \pi )^{2k-1}   \Lambda(\G_{2k};2k-1) =  -\frac{(2k-2)!}{2}  \,  \zeta(2k-1)  $$
 and  $\Lambda(\G_{2k};2k-1) = (-1)^k \Lambda(\G_{2k};1)$. 
 This confirms conjecture \ref{conjEisasMZV}  (and theorem \ref{thmSaad}) in the case of length one. The numbers  $\zeta(2k-1)$ are  periods of a simple extension $$  0 \To \Q  \To M_{\G_{2k}}  \To \Q(1-2k) \To 0 $$
 of Tate motives. This  extension (or a Tate twist of its dual), can be realised as a subquotient of the affine ring of the relative completion $\Or(\GG_{1,1})$   and  corresponds to the Eisenstein series (it is dual to $\be_{2k} \otimes V_{2k-2}^{dR} \in \uu_{1,1}^{dR}$).   This extension is the source of the `zeta elements' $\sigma^{dR}_{2k-1}$.
 
   \subsection{General case} 
   For any   modular form $f$ of weight $w$ set:
  \begin{equation} \label{Pfdefin} P_f(y) = \sum_{k=1}^{w-1} i^{w-k-1} \binom{w-2}{k-1}  \Lambda(f ;k)  y^{k-1} \ .
  \end{equation}
     The function $P_f$ can be interpreted as the value of a canonical cocycle  $C_f$  in a certain cochain complex. The cocycle relations imply functional relations for $P_f$.  
         When $f$ is cuspidal,  these  are equivalent to the   period polynomial equations:
   \begin{eqnarray} P_f(y)+ y^{w-2} P_f(-y^{-1}) & = & 0  \nonumber \\
      P_f(y) + (1-y)^{w-2}P_f \left( \frac{1}{1-y} \right) + y^{w-2} P_f \left(\frac{y-1}{y}\right) &=&  0  \ .  \nonumber 
      \end{eqnarray}
   The first equation follows from the functional equation of $\Lambda$.  A variant  of these equations is satisfied for an Eisenstein series  $f=\G_{w}$  (\cite{MMV}, \S7).   Kohnen and Zagier's  `extra relation'   \cite{KohnenZagier}, which expresses orthogonality of cusp forms to Eisenstein series,  is a consequence of  transference equations (\cite{MMV}, \S8).  
   Manin showed \cite{Ma0}  that when $f$ is a Hecke eigenform, the function $P_f$ is an eigenfunction for a certain action of Hecke operators. From this he deduced that 
   \begin{equation} \label{ManinrelationforPf}  P_f = \omega_f^+  P_{f,+}  + i\, \omega_f^{-} P_{f,-}
   \end{equation} 
   where $\omega_f^+, \omega_f^- \in \R$ and $P_{f,+}, P_{f,-} \in K_f[y]$ where $K_f$ is the  field generated by the Fourier coefficients of $f$.   
   These polynomials  play a key role in the Ihara-Takao relation and its generalisations.

  \begin{example}
       Let $\Delta$ denote the  unique normalised cusp form of weight 12:
 $$\Delta = q\prod_{n\geq 1}(1-q^n)^{24}= \sum_{n \geq 1} \frac{a_n}{n^s}$$
 where $a_1=1, a_2=-24, a_3= 252, \ldots$.
 Equation \eqref{ManinrelationforPf} holds with  
  \begin{eqnarray}  \label{PDelta}
 P_{\Delta, +}   &= &   \frac{36}{691} (y^{10}-1) + y^2-3  y^4+3 y^6-y^8 \\ 
 P_{\Delta,-} & = & 4y -25y^3+42  y^5-25 y^7+4 y^9 \nonumber 
 \end{eqnarray} 
 where $\omega_{\Delta, +}, \omega_{\Delta,-}$ take the numerical values
$$\omega_{\Delta,+} =0.114379022438848\ldots \quad , \quad    \omega_{\Delta,-} =  0.009269276162370\ldots \ . $$
The numbers $(2\pi i)^{11} \omega_{\Delta, \pm}$ are examples of totally holomorphic multiple modular values.
 \end{example}

 \subsection{Motivic interpretation} \label{sectsubMotInterp}  The Manin relation \eqref{ManinrelationforPf}   can be interpreted as an instance of Deligne's conjecture on critical values of $L$-functions.  Associated to a Hecke eigen cusp form $f$ is a rank two motive $M_f$ with coefficients in  the field $K_f$ generated by the Fourier coefficients of $f$ \cite{Scholl}. It has   Hodge numbers of type  $(0,w_f-1)$ and  $(w_f-1,0)$.  The periods of  $M_f(w_f-1)$ are encapsulated by  a two-by-two matrix:
\[
\begin{blockarray}{ccc}
  &  \eta_{f}  &  \omega_{f}\\
\begin{block}{c(cc)}
 P_{f,+}&  \eta_{f,+} & \omega_{f,+}  \\
   P_{f,-}  & i \eta_{f,-}  &  i \omega_{f,-} \\
\end{block}
\end{blockarray} 
 \]
 The periods of $M_f$ are given by the same matrix with all entries multiplied by $(2\pi i)^{w_f-1}$. Here,   $(2\pi i)^{w_f-1} \omega_f,  (2\pi i)^{w_f-1}\eta_f$ are a $K_f$-basis of   the de Rham realisation $\left(M_f\right)_{dR}$ and $P_{f,\pm}$ correspond to a basis for the Betti realisation    $(M_f)_B$.   Poincar\'e duality states that 
 $$M_f^{\vee} = M_f(w_f-1)$$
 and hence    $\bigwedge^2 M_f = K_f(1-w_f)$. This implies a Legendre-style  relation:
  $$    i \eta_{f,+} \omega_{f,-} -  i  \omega_{f,+}   \eta_{f,-}   \quad \in \quad   (2 \pi i)^{1-w_f}\, K_f^{\times} \ . $$

  The holomorphic class $\omega_f \in F^{w_f-1} $ corresponds to the modular form $f$ itself; the class $\eta_f$ can be generated by 
 a modular form of the `second kind', and 
  is only well-defined up to adding a $K_f$-rational multiple of $\omega_f$.

  Deligne's conjecture predicts that the values of $\Lambda(f,s)$ at a  critical point $s$  (i.e., $s=1,\ldots, w_f-1$) are proportional, depending on the parity of $s$, to either of the two holomorphic periods $\omega_f^+, \omega_f^-$. This is indeed  confirmed by Manin's equation \eqref{ManinrelationforPf}.  
 
 Let us rephrase this in a different manner. The values $\Lambda(f;n)$ are periods of an extension of mixed `motives' of the following form:
 $$ 0 \To K_f(n) \To \mathcal{E}  \To  M_{f}(w_f-1)\To 0\ .$$
 It (or rather, it tensored with the  representation $V^{dR}_{w_f-2}$   of $\SL_2$ of dimension $w_f-2$ ) arises as a subquotient of $\Or(\GG_{1,1})$ and  corresponds to the form $f$ (dual to the class $\be'_f \otimes V_{w_f-2}^{dR}$).  
 However,  unlike the case of the Eisenstein series, this extension necessarily  splits when  $1 \leq n \leq w_f-1$.
 This is why  $\Lambda(f;n)$, for $1\leq n\leq w_f-1$, is  in fact a period of  the summand  $M_f(w_f-1)$.

 \subsection{Quasi-periods and non-totally holomorphic MMV's} \label{sectQP}
  In order to obtain the two other non-classical periods $\eta_{f,+}, i \eta_{f,-}$ (which could be called `quasi-periods' following the terminology for algebraic curves), one needs to consider modular forms of the second kind  with poles, for example,  at the cusp. 
  
  Consider  the Ramanujan cusp form $\Delta$. There exists a unique weakly holomorphic  modular form $\Delta'$ of weight $12$ which has a simple pole at the cusp, and whose Fourier coefficients $a_0$, $a_1$ vanish. It is a weak Hecke eigenform with the same eigenvalues as $\Delta$ and has integer coefficients.  Explicitly,  
 $$ \Delta'= q^{-1} +47709536\, q^2+39862705122\, q^3+7552626810624\, q^4+ \ldots$$
Its periods  are  \cite{CNHMF3}
 $$\eta_{\Delta, +}  =   211.113366616704346\ldots \quad , \quad  \eta_{\Delta, -} = 17.055972753974248\ldots\ ,$$
and  one checks  that $    i \left(\eta_{\Delta, +} \omega_{\Delta, -} - \omega_{\Delta,+}  \eta_{\Delta,-}\right)  = 10!  \times (2 \pi i)^{-11}$.

 In general, the periods $ (2\pi i)^{w_f-1} \eta_{f,\pm}$ are  examples of non totally-holomorphic multiple modular values.   Since they are not canonically defined, they
 fall  outside the remit of the conjectures of special values of $L$-functions.  
 We shall not say any more about non-totally holomorphic multiple modular values, as very few are known. In principle, one should  be able to compute them using  the methods of Luo \cite{Luo}.

\section{Examples of double Eisenstein integrals} \label{secMMVex2}

We first discuss concrete examples of multiple modular values in length two before explaining their general structure in the next section. 
\subsection{Examples in low weights} 
These  only involve multiple zeta values.

\begin{example} (Two Eisenstein series of weight $4$).   
Write  $\Lambda_{i,j}= \Lambda(\G_4,\G_4;i,j)$ for brevity,  where $1\leq i,j\leq 3$. 
One can prove that: 
$$\Lambda_{1,1}=  \frac{\zeta(3)^2}{2^7 \,\pi ^6}  \quad , \quad   \Lambda_{1,2}=\frac{\zeta(3)}{ 2^5 3^2 5 \, \pi^3} - \frac{5\, \zeta(5)}{2^8 3 \,\pi^5}  \quad , \quad   \Lambda_{1,3}=  \frac{\zeta(3)^2}{2^6  \pi^6}-\frac{11. 19}{2^{10} 3^{4}  5^{3}}$$

$$\Lambda_{2,1}= -\frac{\zeta(3)}{2^83.5 \,\pi^3}+\frac{5 \, \zeta(5)}{2^83 \, \pi^5}  \quad , \quad \Lambda_{2,2}= \frac{1}{2^{11} 3^4}  \quad , \quad  \Lambda_{2,3} =\frac{\zeta(3)}{2^5 3^25\, \pi^3} -\frac{ 5\, \zeta(5)}{ 2^8 3 \, \pi^5} $$ 

$$\Lambda_{3,1} = \frac{11.19}{2^{10} 3^45^3}     \quad , \quad  \Lambda_{3,2}= - \frac{\zeta(3)}{2^8 3. 5 \,\pi^3}+ \frac{5 \, \zeta(5)}{2^8 3\, \pi^5}   \quad , \quad \Lambda_{3,3} = \frac{\zeta(3)^2}{ 2^7 \pi^6}  $$
In particular,  one can check the symmetry $\Lambda_{i,j} =   \Lambda_{4-j, 4-i}$  and  the shuffle  product relation
$\Lambda_{i,j} +\Lambda_{j,i}   =  \Lambda(\G_4;i)\Lambda(\G_4; j)$.

\end{example}

\begin{example} (Double integrals of two Eisenstein series of weights $6$ and $4$.)  Consider the values of $\Lambda(\G_6,\G_4; n_1, n_2)$ for $1\leq n_1\leq 5$ and $1\leq n_2\leq 3$.   There is a symmetry 
$$\Lambda(\G_6,\G_4; n_1, n_2) =  \Lambda(\G_4, \G_6; 4-n_2, 6-n_1)$$
and shuffle relation 
$$\Lambda(\G_6, \G_4; n_1, n_2) + \Lambda(\G_4, \G_6;n_2,n_1) = \Lambda(\G_6;n_1)\Lambda(\G_4; n_2).  $$
The  corners 
$$\Lambda(\G_6,\G_4;1,1) = \frac{ 3 \,\zeta(3, 5)}{2^6 5 \pi^8}- \frac{503}{2^{11} 3^4 5^3 7}$$
and 
$$\Lambda(\G_6,\G_4;5,3) =   \frac{ 3 \,\zeta(3, 5) }{ 2^6 5  \pi^8} +  \frac{3\, \zeta(3) \zeta(5)}{2^6 \pi^8}  -\frac{503}{ 2^{11} 3^4 5^3 7 }$$  
 give rise to the first non-trivial multiple zeta values in weight $8$, namely $\zeta(3,5)$, which is  conjecturally algebraically independent from the values of the Riemann zeta function at integers $\zeta(n)$.  The remaining values, denoted  by  $\Lambda_{i,j}= \Lambda(\G_6,\G_4;i,j)$,  are:

$$\Lambda_{1,2} = -\frac{\zeta(3)}{ 2^6 3^2 5. 7 \, \pi^3}- \frac{\zeta(5)}{2^8 3 \, \pi^5}+\frac{7 \, \zeta(7)}{2^9\, \pi^7} \quad , \quad  \Lambda_{1,3} = \frac{\zeta(3)^2}{2^6 7  \, \pi^6}+\frac{19.23}{2^{10}3^5 5^27^2}- \frac{3\,  \zeta(3)\zeta(5)}{2^6  \,\pi^8}   $$

$$\Lambda_{2,1} = \frac{\zeta(3)}{2^8 3^2 5.7 \, \pi^3}- \frac{7 \, \zeta(7)}{2^{10} \, \pi^7}  \ , \quad        \Lambda_{2,2}= - \frac{\zeta(3)^2}{2^8 7\, \pi^6}  -\frac{107}{2^{10}3^5 5^2 7^2 } \ ,   \quad  \Lambda_{2,3} = \frac{\zeta(3)}{2^9 3.5.7\, \pi^3}-\frac{\zeta(5)}{2^73.5\, \pi^5} $$

$$\Lambda_{3, 1} = \frac{\zeta(3)^2}{2^7 3.7 \, \pi^6} -\frac{1187}{2^{10}3^6 5^2 7^2}   \  , \quad  \Lambda_{3, 2}= -\frac{\zeta(3)}{2^73^3 7 \, \pi^3}+\frac{\zeta(5)}{ 2^8 3. 5 \, \pi^5} 
 \   , \quad  \Lambda_{3, 3}= \frac{ \zeta(3)^2}{2^7 3.7 \pi^6} -  \frac{1187}{2^{10} 3^6 5^2 7^2} $$

$$\Lambda_{4, 1}= \frac{\zeta(3)}{2^8 3^2 7 \, \pi^3}- \frac{\zeta(5)}{ 2^7 3. 5\, \pi^5} \   , \quad \Lambda_{4, 2}= -\frac{\zeta(3)^2}{2^8 7  \, \pi^6}+\frac{521}{2^{11}3^5 5^2 7^2} \  , \quad \Lambda_{4, 3}=  \frac{\zeta(3)}{2^9 5. 7 \, \pi^3}-\frac{7\, \zeta(7)}{2^{10}\,  \pi^7}$$

$$\Lambda_{5, 1}= \frac{\zeta(3)^2}{2^6 7 \pi^6}+\frac{ 19.23}{2^{10}3^5 5^2 7^2} \quad , \quad \Lambda_{5, 2}= -\frac{\zeta(3)}{2^6 3^2 5. 7\, \pi^3}+\frac{7 \, \zeta(7)}{2^9 \, \pi^7}$$
\end{example}
Continuing in this manner, one can show that all  $\Lambda(\G_a,\G_b; n_1,n_2)$ for $0\leq n_1 <a, 0< n_2<b,$  and $a+b\leq 10$  are multiple zeta values.

\subsection{Examples in weight $12$} The first modular periods appear. 
\begin{example}   \label{Example410} Consider the double Eisenstein integrals $\Lambda(\G_4, \G_{10}; n_1,n_2)$.  There are too many periods to write down in full, but  we find exactly two new periods which are of modular type and not expected to be multiple zeta values, namely:
$$ \pi^{-1}   \Lambda(\Delta; 12)     =   600 \, \Lambda(\G_4,\G_{10};2,5) +   480 \,       \Lambda(\G_4,\G_{10};3,4)   $$
which is a \emph{non-critical}  value of the completed $L$-function of $\Delta$.    This can be proved using the Rankin-Selberg method \cite{MMV} \S9. There is a  new number:
$$c(\Delta;12) = 70 \,  \Lambda(\G_4,\G_{10};3,5) =0.000225126548190262999168981015\ldots \ $$
which is  only well-defined modulo $\Q$.  For example, we might equally well have taken the following quantity as our representative for it
$$\Lambda(\G_4, \G_{10};2,6) = \frac{13}{2^{12}3^5 5^2 7. 11}+\Lambda(\G_4,\G_{10};3,5)$$
which differs by a rational number. See  \S \ref{sectInterp} for an interpretation of these numbers. We now turn to the most interesting multiple zeta values which occur. 
\vspace{0.05in}

The multiple modular value in the  corner satisfies
\begin{equation} \label{f93period} \Lambda(\G_4,\G_{10};1,1)  \overset{?}{=}   \frac{2^2 3^2}{691} c(\Delta;12)  - \frac{3^2 5.7}{2^6} \frac{\zeta_{3,9}}{\pi^{12}} 
\end{equation}
where $\zeta_{3,9}$ is the multiple zeta value of weight $12$ and depth $4$ given by \eqref{f93formula}. Since $c(\Delta;12)$ is of coradical filtration one (see below), this is indeed consistent with  conjecture \ref{conjEisasMZV}, and the stronger versions of it discussed immediately after.

\begin{rem}
Another  interesting value  is
\begin{equation} \label{f37period} \Lambda(\G_4,\G_{10}; 3,1 ) \overset{?}{=}  \frac{4027}{2^8 3^5 5^2 7.11^2} +\frac{c(\Delta;12)}{3^2 5} -\frac{5^2 7}{2^5 11} \frac{\zeta_{3,7} }{\pi^{10}}
     \end{equation}
     where $\zeta_{3,7}$ is  the  multiple zeta value \eqref{zeta73etal}.
The  appearance of this number is explained by `transference' from a different double Eisenstein integral (see  \S\ref{sectTransference}). 
\end{rem}

In the previous equations, $ \overset{?}{=} $ denotes an equality which is true to hundreds of digits, but is presently unproved.  A potential strategy of proof is sketched below. 
\end{example}

\begin{example} \label{Example68}
Consider the double Eisenstein integrals   $\Lambda(\G_6,\G_8;n_1,n_2)$. The values for $0<n_1<6, 0<n_2<8$ satisfy similar properties to the previous example. They  are conjecturally linear combinations of $c(\Delta;12)$ and $\pi^{-1} \Lambda(\Delta;12)$ and multiple zeta values. Conjecture \ref{conjEisasMZV} concerns
\begin{equation}\label{f75period}
\Lambda(\G_6, \G_8; 1,1) \overset{?}{=}  -\frac{2. 3^4}{7.691} c(\Delta;12)  - \frac{3^3 5 }{ 2^7}  \frac{\zeta_{5,7}}{\pi^{12}}
\end{equation} 
where  $\zeta_{5,7}$   can be deduced from \eqref{f93formula} and \eqref{f75andf93}.

\end{example}

\subsection{Interpretation of the modular  periods arising in weight 12}  \label{sectInterp}
Consider  an extension of the following form, where $M_{\Delta}$ was considered in \S\ref{sectsubMotInterp}:
\begin{equation} \label{MDeltaext} 0 \To M_{\Delta} \To \mathcal{E} \To \Q(-12) \To 0\ . 
\end{equation} 
Beilinson's conjecture suggests that such an extension is essentially unique: i.e., there is a one-dimensional space of such extension classes which can arise from the cohomology of  algebraic varieties over $\Q$. We find several equivalent such extensions (tensored with $V^{dR}_{10}$) inside  $\Or(\GG_{1,1})$   and the cocycle relations indeed confirm that they all have the same periods, in accordance with Beilinson's conjecture.  In detail, one  finds that 
$$\mathcal{E}_{dR} = (M_{\Delta})_{dR} \oplus \Q(-12)_{dR}$$
is canonically split by the  Hodge and weight filtrations: it has Hodge types $(0,11), (11,0)$ and $(12,12)$.  The iterated integrals of two Eisenstein series $\G_4, \G_{10}$ correspond to classes in $F^{12} \mathcal{E}_{dR}$. They are dual to   $[\be_4 V^{dR}_2, \be_{10} V^{dR}_{8}] \subseteq \uu^{dR}_{1,1}$. 
Computing the  iterated integrals yields the following period matrix for $\mathcal{E}(11)$: 
 \[
\begin{blockarray}{ccc}
 \begin{block}{(ccc)}
   \eta_+ &   \fbox{$\omega_+$}     &    \fbox{$ \Lambda(\Delta;12)$}  \\ [6pt]
  i \eta_{-}  &  \fbox{$ i \omega_{-}$}  &      90 \, (2 \pi i)  c(\Delta; 12)   \\  [6pt]
  0 & 0 &  \fbox{$2 \pi i$}  \\
\end{block}
\end{blockarray}
 \]

The space of real Frobenius invariants in $\mathcal{E}_B(11)$ is one-dimensional (rows with real entries), but the space of real Frobenius anti-invariants (imaginary rows) is two-dimensional. Changing basis in  $\mathcal{E}_B(11)^{\vee}$ in a way that respects  Frobenius  can  add the third row to the second,   which effectively  modifies  $c(\Delta;12)$ by addition of a rational.  The boxed elements are  well-defined up to scalar multiple. They are all predicted by Deligne's conjecture on special values of $L$-functions, or Beilinson's in the case of the non-critical value $\Lambda(\Delta;12)$.  The three  other non-zero entries are not canonically defined and are not predicted by the standard conjectures on $L$-functions.

Since $\Lambda(\Delta;12) \neq 0$, the  extension \eqref{MDeltaext} is non-trivial, which proves the existence of the  first modular element $\sigma^{dR}_{\Delta(12)}$. The subspace $(M_{\Delta})_{dR}$ is dual to $\be_{\Delta}V^{dR}_{10} \in \uu_{1,1}^{dR}$  (where $\be_{\Delta}$ denotes the pair $(\be'_{\Delta}, \be''_{\Delta})$ ). Therefore $\sigma^{dR}_{\Delta(12)}$ is a derivation which sends 
$$\be_{\Delta} V^{dR}_{10} \quad \mapsto \quad  [\be_4 , \be_{10}] V^{dR}_{10} + \ldots $$
where the right-hand $V^{dR}_{10}$ is a factor of  $V^{dR}_2 \otimes V^{dR}_{8} \cong V^{dR}_{10} \oplus V^{dR}_8 \oplus V^{dR}_6$.

Equivalent extensions show up in many other places inside $\Or(\GG_{1,1})$.

\subsection{Modular incarnation of the Ihara-Takao  relation}  \label{sectModularIncarnation}
The equation  \eqref{ModularIT} can be viewed  on the level of periods as the identity
$$  9\, \Lambda( \G_4, \G_{10};1,1) +14\, \Lambda(\G_6,\G_8;1,1)   =  -\frac{3^3 5.7}{2^6}\left( \zeta_{5,7}+ 3\, \zeta_{3,9}\right)   \ , $$
as deduced by taking  $9\times$equation \eqref{f93period} and adding  $14 \times$equation  \eqref{f75period}. The modular period $c(\Delta;12)$ cancels out and the expression  in brackets is dual to the terms in the Ihara-Takao identity \eqref{Ihara-Takao}.    Here we see in a very concrete manner the exchange of information between the depth filtration on multiple zeta values, the appearance of modular periods as double Eisenstein integrals, and the Ihara-Takao relation.

\section{Structure of double integrals} \label{sectStructuredouble}

In order to continue our exploration of multiple modular values of length two, we need to express the answers more compactly using generating functions.
\subsection{Generating functions for length two} For $f,g$ modular forms of level one and weights $w_f,w_g$, define a polynomial in two variables by:
$$P_{f,g}(y_1,y_2) = \sum_{k,\ell} i^{w_f+w_g-k-\ell-2}  \binom{w_f-2}{k-1}\binom{w_g-2}{\ell-1}\Lambda(f, g; k ,\ell) y_1^{k-1}y_2^{\ell-1}\ .$$ 
It has bidegree $(w_f-2,w_g-2)$. The shuffle product implies that
\begin{equation}\label{Pfgshuffle} P_{f,g}(y_1,y_2) + P_{g,f}(y_2,y_1) =
P_{f}(y_1)   P_{g}(y_2)
\end{equation}
where the terms on the right-hand side are  \eqref{Pfdefin}.  
 In order to unpackage $P_{f,g}(y_1,y_2)$ into simpler pieces  we must use a little representation theory of $\mathrm{SL}_2$.  First define
$$\widetilde{P}_{f,g}(X_1,X_2,Y_1,Y_2) =   X_1^{w_f-2} X_2^{w_g-2} P_{f,g} \left( \frac{Y_1}{X_1}, \frac{Y_2}{X_2}\right)  $$
to be  the homogeneous version in four variables. 
 One shows that this polynomial is a value of a canonical  $\mathrm{SL}_2(\Z)$-cochain $C_{f,g}$  which satisfies the equation  \cite{MMV}, \S5 
$$\delta C_{f,g} = C_f \cup C_g\ .$$
This, together with some known information  about the value of cocycles at infinity   implies many (but not all) relations between the coefficients $P_{f,g}$. We shall not say any more about this here and turn straight to examples. To this end, consider
$$\partial =  \frac{\partial}{\partial X_1}\frac{\partial}{\partial Y_2} -\frac{\partial}{\partial Y_1}\frac{\partial}{\partial X_2}$$
which is an $\mathrm{SL}_2$-equivariant differential operator  of degree $-2$. 
The two-variable polynomial $P_{f,g}$  is uniquely determined by a finite   sequence of  simpler polynomials in a single variable of degree $w_f+w_g-2k-4$ defined for $k\geq 0$ by:
\begin{equation} 
\delta^k P_{f,g}  =    \left(\partial^k  \widetilde{P}_{f,g}\right)(1,1,y,y)\ .\end{equation} 
The  polynomial $\delta^0 P_{f,g}$ is simply the diagonal $P_{f,g}(y,y)$.  The polynomials $\delta^k P_{f,g}$  vanish for all  $k\geq \min\{w_f-2,w_g-2\}$.

\begin{example}  The discussion of example  \ref{Example410}  can be summarised using the generating function  $P_{\G_4,\G_{10}}$. It  has the following shape
$$ P_{\G_4,\G_{10}} = \frac{ \zeta_{3,9} }{\pi^{12}} c_0 \, B^{(0)} + \frac{\zeta_{3,7}}{\pi^{10}}c_1\, B^{(2)}   +   
  \frac{\Lambda(\Delta;12)}{\pi}c_2  \,P^{(0)}_{\Delta,-}  + c(\Delta;12)c_3 \,P^{(0)}_{\Delta,+} + \, R
$$
where $c_0,c_1,c_2,c_3$ are known rational numbers, $B^{(0)},B^{(2)}$ are  the unique solutions to
$$\begin{array}{lcl}
  \delta^0 B^{(0)}  =  y^{10}-1  & \qquad   &    \delta^0 B^{(2)} = 0 \\
  \delta^1 B^{(0)} =  0 &   &   \delta^1 B^{(2)} =  0   \\
 \delta^2 B^{(0)} = 0  &   &  \delta^2 B^{(2)} =   y^6-1  
\end{array}
$$ ($B^{(0)}= y_1^2y_2^8-1$ and  $B^{(2)} = (y_1-y_2)^2(y_2^6-1)$), 
and $P^{(0)}_{\Delta,\pm}$ are  the unique solutions to
$$\begin{array}{lcl}
 \delta^0  P^{(0)}_{\Delta,-} = P_{\Delta, -} & \qquad   &    \delta^0  P^{(0)}_{\Delta,+} = P_{\Delta, +} \\
\delta^1  P^{(0)}_{\Delta, -} =  0 &   &  \delta^1  P^{(0)}_{\Delta,+} =  0   \\
 \delta^2  P^{(0)}_{\Delta,-} = 0  &   &   \delta^2  P^{(0)}_{\Delta,+} =   0  
\end{array}
$$
where $P_{\Delta, \pm}$ are the odd and even period polynomials  \eqref{PDelta}   associated to $\Delta$.
The terms $B^{(-)}$ can be interpreted as coboundaries  in a cochain complex.  
The remainder term $R$ involves only single zeta values and products. It is mostly determined from  $P_{\G_4}$ and $P_{\G_{10}}$ via \eqref{Pfgshuffle},
although some new single zeta values can occur, including  $\frac{\zeta(11)}{\pi^{11}}$. \end{example}

\begin{example} Example \ref{Example68} can be reformulated in a similar way.
 It is:
 $$ P_{\G_6,\G_{8}} = \frac{ \zeta_{5,7} }{\pi^{12}} c_0 \, B^{(0)} + \frac{\zeta_{5,3}}{\pi^{8}}c_1\, B^{(4)}  +   
  \frac{\Lambda(\Delta;12)}{\pi}c_2  \,   P^{(0)}_{\Delta,-}    + c(\Delta;12)c_3 \,  P^{(0)}_{\Delta,+}        + \, R
$$
where $c_0,c_1,c_2,c_3\in \Q$, the $B$ are the unique polynomials of bidegree $(4,6)$ satisfying
$$ \delta^i B^{(k)} = \begin{cases}  y^{10-2i}  -1   \qquad \hbox{ if } i =k    \\
0  \qquad  \qquad  \qquad \hbox{ if } i \neq k   \end{cases}
$$
 (e.g, $B^{(0)}= y_1^4y_2^6-1$), 
and $P^{(0)}_{\Delta,\pm}$ are the unique polynomials of bidegree $(4,6)$ with
$$\delta^i P_{\Delta,\pm}^{(0)} = \begin{cases} P_{\Delta, \pm}  \qquad \hbox{ if } i =0    \\
0    \qquad  \qquad \hbox{ if } i >0   \end{cases} $$
and 
$R$ involves only single zeta values and products. 

\end{example}

\subsection{General case in length two}
Let $f,g$ be as above. The generating function $P_{f,g} \in \C[y_1,y_2]$   has bidegree $(w_f-2,w_g-2)$ and takes  the form
$$P_{f,g} = \sum_{k\geq 0}  \left(  c_{f,g;\G_{d_k}} \, B^{(k)}  +  \sum_{h}  c_{f,g;h}^+ \, P^{(k)}_{h,+} + \sum_{h} i\, c_{f,g;h}^-\,  P^{(k)}_{h,-} \right)  +  R$$
where $d_k= w_f+w_g -2 -2k$, the second  and third sums range over a basis $h$ of  Hecke eigen cusp forms $h$ of weight $d_k$, and  $P_{h,\pm}$ are (a choice of) associated
odd and even period polynomials \eqref{ManinrelationforPf} with coefficients in the field $K_f$ generated by the eigenvalues of $h$.
  The coboundary polynomials $B$ are the unique solutions to the equations:
$$
\delta^j  B^{(k)} \quad =  \quad  
\begin{cases}  
0  \qquad \qquad \quad  \hbox{ if  } \qquad  j\neq k  \\
 y^{d_k-2}  -1  \quad  \ \hbox{ if } \qquad  j=k 
\end{cases} \ .
$$
The polynomials  $ P^{(k)}_{h,\pm}$ are the unique solutions to the equations
$$
\delta^j  P^{(k)}_{h,\pm} \quad =  \quad  
\begin{cases}  
0  \qquad \qquad \hbox{ if  } \qquad  j\neq k  \\
 P_{h,\pm} \   \qquad \hbox{ if } \qquad  j=k
\end{cases} \ .
$$
 The remainder $R$ is completely (but not uniquely!) determined  from  $P_f$ and $P_g$. Its coefficients are linear combinations of  products of the coefficients of the latter. 
 
 All the interesting information is contained in the numbers  
 $$c_{f,g;\G_{d_k}} \ , \   c_{f,g;h}^+ \ , \   c_{f,g;h}^-\  \qquad \in  \quad \R$$
 which are only well-defined modulo products and periods of lower unipotency filtration.  The reason for the notation 
 $c_{f,g;\G_{d_k}}$ is that the cocycle of  the Eisenstein series $\G_{d_k}$ is  Poincar\'e dual to the coboundary cocycle $\delta^k \, B^{(k)} = y^{d_k-2} -1$. 

 \begin{rem} The terms  in the above   formula for $P_{f,g}$ are a Hecke eigenbasis for the group cohomology 
$H^1(\SL_2(\Z),    V^{dR}_{w_f-2} \otimes V^{dR}_{w_g-2} ) $.
 \end{rem} 

\begin{example} (Modular Ihara-Takao relation on generating functions). 
The polynomials $\delta^0 P_{\G_4,\G_{10}}$ and $\delta^0 P_{\G_6,\G_8}$ 
 both have  degree  $10$. All modular periods $\Lambda(\Delta;12), c(\Delta;12)$  cancel out in the  linear combination  $ 9\, \delta^0 P_{\G_4,\G_{10}}  +  14\, \delta^0 P_{\G_6,\G_8}$ (compare \S\ref{sectModularIncarnation}).  \end{example} 

\subsection{Transference principle in length two}  \label{sectTransference}
 In length one, transference  is equivalent to the well-known orthogonality  of cocycles of Hecke eigenforms.  In  length two, transference is a relation  between linear combinations of multiple modular values
$$\Lambda(f;a_1) \Lambda(g,h;a_2,a_3)  \quad \hbox{ and }  \quad \Lambda(f,g;b_1,b_2) \Lambda(h;b_3)$$
where $f,g, h$ are fixed Eisenstein series or cusp forms and $a_1,a_2,a_3,b_1,b_2,b_3$ can vary within the allowed range. 
Heuristically,  we obtain relations of the form:
\begin{eqnarray}
 c_{\G_a, \G_b; \G_c} &  \sim &  c_{\G_b, \G_c; \G_a}  \nonumber \\
c_{f,m;\G_a} & \sim & \left(  c^+_{m,\G_a; f} \,\omega_{f,-} + c^-_{m,\G_a; f} \, \omega_{f,+}\right)
 \nonumber \\
\left(  c^+_{f,m;g}\, \omega_g^{-}  +   c^-_{f,m;g}\, \omega_g^{+} \right)  & \sim & \left(  c^+_{g,m; f} \,\omega_{f,-} + c^-_{g, m; f} \, \omega_{f,+}\right)
 \nonumber 
\end{eqnarray} 
where $\G_a,\G_b,\G_c$ are Eisenstein series, $f,g$ are cuspidal eigenforms and $m$ is a modular form and $\sim$ denotes  a  relation modulo products and periods of lower coradical filtration. These equations also hold for non-totally holomorphic multiple modular values.

\begin{example} 
The first transference equation between the three Eisenstein series $\G_4, \G_8, \G_{10}$ is an explicit relation between the multiple modular values of $\Lambda(\G_4, \G_{10})$ and those of $\Lambda(\G_4, \G_8)$. 
For example, one can deduce the identity:
$$14\,  \Lambda(\G_4,\G_{10};3,5) - 9\, \Lambda(\G_4,\G_{10};3,1)  =   \frac{2^25.7}{11} \Lambda(\G_4,\G_{8};1,1) + \frac{157}{2^7 3^3 5^2 11^2} \ .  $$
Now conjecture \ref{conjEisasMZV} predicts that 
$$ \Lambda(\G_4, \G_8;1,1)  \overset{?}{=}   \frac{3^2 5}{2^7}\frac{\zeta_{3,7}}{\pi^{10}}  - \frac{83}{2^{10} 3^2 5. 7^2 11}$$
and so we can understand the appearance \eqref{f37period} of the multiple zeta value $\zeta_{3,7}$ in $\Lambda(\G_4, \G_{10};3,1)$  as having been transferred from $\Lambda(\G_4, \G_8;,1,1)$.   Conjecture \ref{conjEisasMZV} and  transference  are sufficient to explain all appearances of multiple zeta values (as opposed to products of simple zeta values) amongst double Eisenstein integrals. 
\end{example}

\begin{example} Transference  between $\G_4, \G_{10}$ and $\Delta$ implies, amongst other things,  that the quantity 
$ c^+_{\G_4,\G_{10}; \Delta} \,\omega_{\Delta,-} + c^-_{\G_4,\G_{10}; \Delta} \, \omega_{\Delta, +}$
which is a linear combination of multiple modular values  of $\Lambda(\Delta)\Lambda(\G_4,\G_{10})$,
 is  transferred to $c_{\Delta, \G_4; \G_{10}}$, which is a linear combination of multiple modular values  of $\Lambda(\Delta,\G_4)\Lambda(\G_{10})$.  
For example:
\begin{multline} (2^3 3^3 5^2 11)^{-1} \left(   44\,  \Lambda(\Delta,\G_4;11,1) -110 \, \Lambda(\Delta,\G_4;10,2) 
+ 15\, \Lambda(\Delta, \G_4; 6, 2)   \right)\quad  =  \nonumber   \\   
 \Lambda(\Delta;3) \left( 30   \Lambda(\G_4,\G_{10};2,5) +24 \, \Lambda(\G_4;\G_{10};3, 4) \right) -35 \Lambda(\Delta;2) \Lambda(\G_4, \G_{10};3,5) 
\end{multline}
where the  right hand side is  
$$ \left(\frac{\omega_{\Delta, +} }{900 \, \pi }  \Lambda(\Delta;12) - \frac{\omega_{\Delta, -}}{5}  c(\Delta;12)\right) $$
Transference explains why periods of the cusp form $\Delta$ can occur as double Eisenstein integrals in weight $12$, without appealing to the Rankin-Selberg method. 
\end{example}

\subsection{Summary of possible extensions in length two}
Length two multiple modular values can give rise to the following `submotives'  (up to Tate twists) in the affine ring $\Or(\GG_{1,1})$ of the relative completion of the fundamental group of $\mathcal{M}_{1,1}$:

\begin{enumerate}
\item (Three Eisenstein series). The numbers $c_{\G_a,\G_b;\G_c}$ are periods of a biextension of Tate objects. We therefore expect them  to be multiple zeta values of the form $\zeta_{2k+1,2\ell+1}$ (definition \ref{defnzetaprescribed}). This is the content of  conjecture \ref{conjEisasMZV}.
\vspace{0.05in}
\item (One cusp form $f$, and two Eisenstein series).  The numbers $c_{f,\G_a;\G_b}$ or $c_{\G_a,\G_b;f}$ are periods of extensions of the form  we have already discussed:
$$0 \To M_f(d) \To \mathcal{E} \To \Q \To 0 $$
where $d$ is non-critical. One period of this extension is a non-critical $L$-value $\Lambda(f;d)$, another is a `non-standard' period of the form $c(f;d)$. 
\vspace{0.05in}
\item (Two cusp forms $f,g$, and one Eisenstein series).  The numbers $c_{f,g;\G_a}$ or $c_{f, \G_a;g}$ are periods of Rankin-Selberg extensions of the form
$$0 \To M_f \otimes M_g(d) \To \mathcal{E} \To \Q \To 0 $$
where $d$ is `semi-critical', i.e., where the rank of the extension group is $1$.  One of its periods is $\Lambda(f\otimes g; d)$, but there are 3 other non-standard periods. 
\vspace{0.05in}
\item (Three cusp forms $f,g,h $).  The numbers $c_{f,g;h}$  are periods of triple Rankin-Selberg extensions of the form
$$0 \To  M_f \otimes M_g\otimes M_h(d) \To \mathcal{E} \To \Q \To 0 $$
where $d$ corresponds to a central critical $L$-value, i.e.,  $\M_f \otimes M_g\otimes M_h(d)$ has weight one, corresponding to the exceptional case in Beilinson's conjecture.
\end{enumerate}

 \emph{Acknowledgements}. The author is  partially supported by ERC grant  724638.  Many thanks to M. Kaneko, M. Luo, and A.  Saad  for corrections. 
\\

\noindent  Francis Brown, \\
All Souls College, \\
Oxford, \\
OX1 4AL,\\
United Kingdom.  \\
\texttt{francis.brown@all-souls.ox.ac.uk}

\bibliographystyle{plain}
\bibliography{main}

\end{document}